\documentclass[12pt, a4paper]{article}
\usepackage{amsmath}
\usepackage{amssymb}
\usepackage{latexsym}
\usepackage{braket}

\newtheorem{thm}{Theorem}[section]

\newtheorem{prop}[thm]{Proposition}
\newtheorem{lemma}[thm]{Lemma}

\newtheorem{col}[thm]{Corollary}

\def \qed {\hfill \hbox {\rule [-2pt]{3pt}{6pt}}}
\title{A representation theorem on a filtering model with first-passage-type stopping time}
\author{Takenobu NAKASHIMA}
\begin{document}
\label{firstpage}
\maketitle
\begin{abstract}
We present a representation theorem for a filtering model with first-passage-type stopping time. 
The model is constructed from two unobservable processes and one observable process that is under the influence of two unobservable processes.
A filter is constructed using Brownian motion in the observable process and a first-passage-type stopping time in an unobservable process.
Though our theorems are similar to those of Nakagawa\cite{Nakagawa}, we do not use pinned Brownian motion measure, which is difficult to deal with.
In addition, we describe a representation theorem for another filtration that was not discussed by Nakagawa\cite{Nakagawa}.
\end{abstract}
%
\section{Introduction}
Duffie and Lando \cite{DuffieLando} studied the implications of imperfect information for 
the term structures of credit spreads on corporate bonds. They assumed that the bond investor
 could not observe the issuer's assets directly, and could receive only periodic and imperfect accounting 
information. They then derived a relationship between the volatility of the issuer's asset value and
 its hazard rate. 
Their model is a kind of filtering model.
Jeanblanc and Valchev \cite{Jeanblanc} examined three types of information related to a company's unlevered asset value on the secondary
bond market: the classical case of continuous and perfect information,  observations of past and
contemporaneous asset values at selected discrete times, and observations of contemporaneous
asset values at discrete times. 
In their model, although bond holders receive information about contemporaneous and past asset values in the second type of information,
 they receive only contemporaneous information in the third type. 
Jarrow, Protter and Deniz \cite{JarrowProtterDeniz} provided an alternative credit risk model based on information
reduction, whereby the market only observes the company's asset value when
it reaches certain levels, interpreted as changes significant enough for the company's
management to make a public announcement.
Nakagawa\cite{Nakagawa}
constructed a filtering model based on a default risk,
and derived representation formulas under conditions of imperfect information.
He analyzed the properties of processes under $\nu^{u,x_2}_{0,x_1}$, 
which is a probability measure on $C([0,u];\mathbf{R})$, and  the law of 
 Brownian motion $B_t$ conditioned to start from
$x_1>0$, stay in $(0,\infty)$ for $s\le u$ and reach $x_2>0$ at time $u$ under $P$.
However, because this measure is difficult to deal with, we
present representation formulas that do not use the measure $\nu$.
In this paper, we refer to the ``first-passage-type stopping time'' instead of a ``default time'', because our focus is solely on the mathematical perspective of a filtering model.

First, we present a representation theorem for a filtration with first-passage-type stopping time. In this part, we do not use a filtration model.

Let $(\Omega ,{\cal B},P,\{ \mathcal{B}_t \}_{t \ge 0})$ be a complete filtrated probability space,
and assume that the filtration $\{ \mathcal{B}_t \}_{t\ge 0}$ satisfies the usual conditions.
Let $B_t$, $\hat{B}_t$ and $W_t$ be independent
$\mathcal{B}_t $-Brownian motions with values in $\mathbf{R}$,$\mathbf{R}^d$ and $\mathbf{R}$
 respectively.
We denote  the right continuous filtration generated by a continuous stochastic process $X$ as $(\mathcal{G}^{X}_t)$.
For example,
 ${\mathcal G}_t^B$ $= \bigcap_{u>t} \sigma\{B_s, s \le u\}$.
Let $a>0$, $B^a_t=a+B_t$, $\tau^a = \inf \{ t>0;\; B^a_t = 0\}$, $N^a_t = 1_{\{\tau^a \le t\}}$ and 
$\mathcal{F}^W_t= \bigcap_{u>t} ( \mathcal{G}^W_u \vee \sigma \{ \tau^a \wedge u \})$.
Let
$q_{a}(t)=\int^{\infty}_t \frac{a}{\sqrt{2\pi s^3}}\exp(-\frac{a^2}{2s})ds$
and
$\lambda_a(t)=-\frac{d}{dt} \log q_{a}(t)=\frac{a}{\sqrt{2\pi t^3}}q_{a}(t)^{-1}\exp(-\frac{a^2}{2t})$.
Then
$P[\tau^a > t ] = q_{a}(t)=e^{-\int^t_0 \lambda_a(u) du}.$
Let $\gamma_a(t)$ be the density of $\tau^a$. Then, we have
\begin{equation}\label{gamma}
\gamma_a(t) dt = P[\tau^a \in dt] = \lambda_a(t) e^{-\int^t_0 \lambda_a(u) du} dt
. 
\end{equation}
We can also see that
$$
M^a_t=N^a_{t}-\int^t_0 (1-N^a_{s})\lambda_a(s) ds
$$
is $\mathcal{F}^W_t$-martingale.

Let $g(t,x)$ and $\Phi(t,x)$ be the density and distribution, respectively, of the Brownian motion $B_t$.
Hence, $g(t,x)$ and $\Phi(t,x)$ can be written as follows.
\begin{equation}\label{g}
g(t,x) = \frac{1}{\sqrt{2\pi t}}\exp ( - \frac{x^2}{2t} ),\;
\Phi(t,x)=\int^x_{-\infty}g(t,y)dy,\; x\ge 0,\; t>0.
\end{equation}
We note that
$$
\frac{\partial g}{\partial x}(t,x) = -\frac{x}{t} g(t,x), 
\;\;\; 
\frac{\partial^2 g}{\partial x^2}(t,x) 
=2\frac{\partial g}{\partial t}(t,x)
= \frac{x^2 -t}{t^2} g(t,x)
.
$$
We denote as $\mathcal{L}^p$, $p \in (1,\infty)$, the space of 
$\{\mathcal{B}_t\}$-progressively measurable functions $\varphi$ such that
 $E[\int^T_0 |\varphi|^p_s ds] < \infty$ for any $T>0$, and write $\mathcal{L}^{p+} = \bigcup_{q>p} \mathcal{L}^q, p \ge 1$. 
For $t>s$, let 
\begin{eqnarray}
H^{(k)}_a(t,s;f) 
=
E[1_{\{ \tau^a >s\} } f_s 
\frac{\partial^k g}{\partial x^k}(t-s,B^a_s)|\mathcal{G}^{W}_{s}], \; f \in \mathcal{L}^{1+},\label{H}\\
\; k=0,1,2, \nonumber \\
\hat{H}^{(k)}_a(t;f) = \int^t_0 H^{(k)}_a(t,u;f)du, \; f \in \mathcal{L}^{\frac{4}{3-k}+},\; k=0,1,2, \label{hatH}\\
\bar{H}_a(t;f)= e^{\int^t_0 \lambda_a(r) dr}
\{\hat{H}^{(2)}_a(t;f)+2 \lambda_a(t) \hat{H}^{(0)}_a(t;f) \}, \label{barH}\\
\; f \in \mathcal{L}^{4+}, \nonumber\\
U_a(t,s;f) = E[1_{\{\tau^a>s\}}f_s (2\Phi(t-s,B^a_s)-1)|\mathcal{G}^W_s],  \; f \in \mathcal{L}^{1+} \label{U}\\
\bar{U}_a(t,s;f) = e^{\int_{0}^{t} \lambda_a(r)dr} 
\{ H^{(1)}_a(t,s;f) + \lambda_a(t) U_a(t,s;f)\}, \label{barU}\\
 \;  f \in \mathcal{L}^{1+} .\nonumber
\end{eqnarray}
We will show that these are well defined in Section \ref{sec:2}.
Thus we have the following theorem.
\begin{thm}\label{th1}
{\rm (1)} 
For any $t,T\ge0$ and $f \in \mathcal{L}^{4+}$,
$$
E[\int^{T}_{0}f_sdB_s|\mathcal{F}^W_{t}]
=
-\int^{t}_{0}\bar{H}_a(s;f 1_{(0,T]}(\cdot))\lambda_a(s)^{-1} dM^a_s.
$$
\\
{\rm (2)} 
For any $t\ge0$ and $f \in \mathcal{L}^{4+}$,
$$
E[\int^t_0 f_sds|\mathcal{F}^W_t]
=
\int^t_0 E[f_s|\mathcal{F}^W_s] ds
-
\int^t_0 \left(\int^s_0 \bar{U}_a(s,r;f) dr\right)
\lambda_a(s)^{-1}
dM^a_s.
$$
\\
{\rm (3)} 
For any $t,T\ge0$ and $f \in \mathcal{L}^{6+}$,
$$
E[\int^T_{0} f_s dW_s|\mathcal{F}^W_{t}]
=
\int^{T \wedge t}_0 E[f_s |\mathcal{F}^W_s]dW_s
-\int^t_0\left(\int^s_0\bar{U}_a(s,r;f1_{[0,T]}(\cdot))dW_r\right)
\lambda_a(s)^{-1} dM^a_s.
$$
\\
{\rm (4)} 
For any $t\ge0$,$\hat{f}_i, \in \mathcal{L}^{2+},   i=1,\cdots,d,$
$$
E[\hat{f}^i_s d\hat{B}^i_s|\mathcal{F}^W_t]=0.
$$
\end{thm}

Second, we consider a representation theorem with a filtering model. 
The quantities $X$, $Z$, and $Y$ are the same as those considered by Nakagawa\cite{Nakagawa},
and are called the main system, sub-system and observation, respectively, in his paper.
Let $X$ and $Z$ be solutions of the following stochastic differential equations under $P$:
\begin{eqnarray*}
dX_t &=& dB_t+b_0(t,X_t,Z_t)dt,\;\;\;\;\;\; X_0=x_0>0,\\
dZ_t &=& \sigma_1(t,X_t,Z_t)d\hat{B}_t+b_1(t,X_t,Z_t)dt,\;\;\;\;\;\;Z_0=z_0 \in \mathbf{R}^{N},
\end{eqnarray*}
where 
$b_0:[0,\infty) \times \mathbf{R} \times \mathbf{R}^{N} \rightarrow \mathbf{R}$,
$\sigma_1:[0,\infty)\times \mathbf{R} \times \mathbf{R}^{N} \rightarrow \mathbf{R}^{N \times d}$
and
$b_1:[0,\infty) \times \mathbf{R} \times \mathbf{R}^{N} \rightarrow \mathbf{R}^{N}$ 
are bounded and continuously differentiable functions.
Let $Y$ be a solution of the   stochastic differential equation,
$$
dY_t = \sigma_2(t,Y_t)dW_t+b_2(t,X_{t\wedge \tau},Y_t)dt,\;\;\;\;\;\;Y_0=y_0 \in \mathbf{R},
$$
where 
$\sigma_2:[0,\infty) \times \mathbf{R} \rightarrow \mathbf{R}$
and
$b_2:[0,\infty)\times \mathbf{R} \times \mathbf{R}\rightarrow \mathbf{R}$
are bounded and continuously differentiable functions.
We assume that there exist some $\epsilon >0$ and  $\sigma_2(t,y)$ satisfying $\sigma_2(t,y) \ge \epsilon$ for any
 $t \in [0,\infty)$, $y \in \mathbf{R}$.
Let
$\tau = \inf \{ t>0;\; X_t = 0\}$,
$N_t = 1_{\{\tau \le t\}}$
and
$\mathcal{F}_t= \bigcap_{u>t} (\mathcal{G}^Y_u \vee \sigma \{ \tau \wedge u \})$.
We now consider changing the probability measure. 
Let 
$\rho_t$ be given by
\begin{eqnarray}\label{rho}
\rho_t &=&  \exp \left(  \int^t_0 b_0(s,X_s,Z_s)dB_s + \int^t_0 \beta (s,X_{s \wedge \tau},Y_s) dW_s \right. \\
&+&\left. \frac{1}{2}\int^t_0 (b_0(s,X_s,Z_s)^2+\beta(s,X_{s \wedge \tau},Y_s)^2)ds
 \right), \nonumber
\end{eqnarray}
where $\beta(t,x,y)=\sigma_2(t,y)^{-1}b_2(t,x,y)$
 and $\widetilde{P}$ is a probability measure on $(\Omega, \mathcal{F})$ given by $d\widetilde{P}=\rho^{-1}_t dP$.
We can see that $\rho, \rho^{-1} \in \bigcap_{p \ge 1} \mathcal{L}^p$ by Novikov's Theorem.
Let
$\widetilde{\rho}_t = \widetilde{E}[\rho_{t} |\mathcal{F}_t]$.
Here, we will denote the expectation under the probability measure $\tilde{P}$ as $\widetilde{E}[\cdot]$.
Let 
\begin{eqnarray*}
\widetilde{B}_t &=& B_t+\int^t_0 b_0(s,X_s,Z_s)ds,\\
\widetilde{W}_t &=& W_t+\int^t_0 \beta(s,X_{s \wedge \tau},Y_s)ds.
\end{eqnarray*}
Then $\widetilde{B}_t$, $\hat{B}_t$ and $\widetilde{W}_t$ are independent $\widetilde{P}$-$\{\mathcal{B}_t\}_{t\in [0,\infty)}$ -Brownian motions.
The stochastic processes $X$, $Z$ and $Y$ are described in the following:
\begin{eqnarray*}
dX_t &=& d\widetilde{B}_t,\\
dZ_t &=& \sigma_1(t,X_t,Z_t)d\hat{B}_t+b_1(t,X_t,Z_t)dt,\\
dY_t &=& \sigma_2(t,Y_t)d\widetilde{W}_t.
\end{eqnarray*}
From the above equations, we can see that $\{\mathcal{G}^X_t\}_{t\in [0,\infty)}$ coincides with the natural filtration
generated by $\{\widetilde{B}_t\}_{t\in [0,\infty)}$.
Because
$d\widetilde{W}_t=\sigma(t,Y_t)^{-1}dY_t$,
we can see that $\mathcal{G}^Y_t=\mathcal{G}^{\widetilde{W}}_t$ and $\mathcal{F}_t= \bigcap_{u>t} (\mathcal{G}^{\widetilde{W}}_u \vee \sigma \{ \tau \wedge u \})$.
In addition, we can see that
$$
\widetilde{M}_t=N^{x_0}_t-\int^t_0 (1-N^{x_0}_{s-})\lambda_{x_0}(s) ds
$$
is $\widetilde{P}$-$\mathcal{F}_t$-martingale.
Let
\begin{equation}\label{I}
I^{(k)}(t,s;f) 
=  \widetilde{E}[1_{\{ \tau >s\} }  
\frac{\partial^k g}{\partial x^k}(t-s,X_s)
\rho_{s-} f_s
|\mathcal{G}^{Y}_{s}],\; t>s, \; k=0,1,2
\end{equation}
for $f \in \mathcal{L}^{2+}$.
Let $\Sigma$ denote the set of $\mathcal{B}$ -adapted continuous processes $F$ for which 
 there exist $f_i$, $i=1,2,3 \in \mathcal{L}^{6+}$ and $(f_4^j),j=1,\cdots,d \in \mathcal{L}^{6+}$such that
\begin{equation}\label{F}
F_t = F_0 + \int^t_0 f_1(s) ds+\int^t_0 f_2(s) dW_s+ \int^t_0 f_3(s) dB_s + \sum_{j=1}^d \int^t_0 f_4^j(s) d\hat{B}^j_s.
\end{equation}\label{operator}

For $F \in \Sigma$, let
\begin{eqnarray}\label{operators}
(\widetilde{D_0} F)_t&=&\beta(t,X_{t \wedge \tau},Y_t)F_{t \wedge \tau}+1_{\{\tau>t\}} f_2(t)  , \\
(\widetilde{D_1} F)_t&=& b_0(t,X_t,Z_t)F_{t\wedge \tau} + 1_{\{\tau>t\}} f_3(t) \nonumber , \\
(\widetilde{D_2} F)_t&=&1_{\{\tau>t\}} f_4(t)\nonumber , \\
(\widetilde{L} F)_t&=& 1_{\{\tau>t\}} f_1(t). \nonumber 
\end{eqnarray}
For $r>s>0$, let
\begin{eqnarray}\label{VandV1}
&\hat{V}(r;F)=& \\
&\widetilde{\rho}_{r-}^{-1} e^{\int^r_0 \lambda_{x_0}(u) du} 
\left(
\hat{V}_1(r,r;F) 
+ \lambda_{x_0}(r) (-(2\Phi(r,x_0)-1)F_0 + \widetilde{E}[1_{\{\tau>r\}}\rho_{r} F_r|\mathcal{G}^Y_r] )
\right)
,& \nonumber\\
&\hat{V}_1(r,s;F) = \nonumber& \\
&\int^s_0 I^{(1)}(r,u;\widetilde{D_0} F) d\widetilde{W}_u
+\int^s_0 \left(I^{(2)}(r,u;\widetilde{D_1} F)+I^{(1)}(r,u;\widetilde{L} F)\right) du
.&\nonumber
\end{eqnarray}
Let
$$
\tilde{\tilde{\lambda}}(s) = \lambda_{x_0}(s)+\hat{V}(s;1),
\;\;\;
\widetilde{\widetilde{M}}_t = N_t - \int^{t}_0 (1-N_s)\tilde{\tilde{\lambda}}(s) ds
$$
and
$$
\widetilde{\widetilde{W}}_t = \widetilde{W}_t - \int^t_0 E[\beta(r,X_{r \wedge \tau},Y_r)|\mathcal{F}_r] dr
.
$$
Then, we will show that $\widetilde{\widetilde{M}}_t$ is $P$-$\mathcal{F}_t$-martingale
 and that $\widetilde{\widetilde{W}}_t$ is a $P$-$\mathcal{F}_t$-Brownian motion. Nakagawa \cite{Nakagawa} also gave $\tilde{\tilde{\lambda}}$
using the measure of a pinned Brownian motion.
We can now state the following representation theorem, which was not given by Nakagawa \cite{Nakagawa}.
\begin{thm}\label{th2}
Let 
$F \in \Sigma$ and $\bar{F}_t=E[F_{t\wedge \tau}|\mathcal{F}_t]$. Then we have the following.\\
{\rm (1)} 
$$
\bar{F}_t
=F_0
+\int^t_{0}  \bar{f}_0(r;F) d\widetilde{\widetilde{M}}_r
+\int^t_{0}  \bar{f}_1(r;F) dr
+\int^t_{0} \bar{f}_2(r;F) d\widetilde{\widetilde{W}}_r,
$$
where
\begin{eqnarray*}
\bar{f}_0(r;F)&=&
-
1_{\{\tau>r\}}
(
\hat{V}(r; F) + \hat{V}(r;1)\bar{F}_{r-}
)
\tilde{\tilde{\lambda}}(r)^{-1}
,\nonumber \\
\bar{f}_1(r;F)&=&
 1_{\{\tau>r\}}
 E[1_{\{\tau>r\}} (\widetilde{L} F)_r|\mathcal{F}_r]
 ,\nonumber \\
\bar{f}_2(r;F) &=&
E[(\widetilde{D_0} F)_r|\mathcal{F}_r]
 - 
E[\beta(r,X_r,Y_r)|\mathcal{F}_r]
\bar{F}_{r-}
. \nonumber
\end{eqnarray*}
\\
{\rm (2)} 
Moreover, if there exist $C>0$ and $\alpha \in (0,1)$  such that $1_{\{|X_t| \le 1\}}1_{\{\tau > t\}}|F_t| \le C|X_t|^{\alpha}$ for $t>0$,
we have
$
\bar{f}_0(r;F)
=
- 1_{\{\tau>r\}}\bar{F}_{r-}
.
$
\end{thm}
The author would like to express his appreciation to Prof. Kusuoka and the referee for their useful suggestions and comments.

\section{Evaluation of integrands}\label{sec:2}
For $f \in \mathcal{L}^1$, $t>s>0$ and $k=0,1,2$, let
\begin{equation}\label{tildeH}
\widetilde{H}^{(k)}_a(t,s;f)=E[1_{\{\tau^a>s\}} |f_s\frac{\partial^k g}{\partial x^k}(t-s,B^a_s)||\mathcal{G}^W_{s}].
\end{equation}
\begin{prop}\label{ineq:HkaVer2_1}
For $q>1$ and $k=0,1,2$, we have
$$
\widetilde{H}^{(k)}_a(t,s;1)
\le
C^{(k)}_1(a)+C^{(k)}_2(q,a)(t-u)^{\frac{-kq-q+2}{2}}
$$
for any $t > u \ge 0$ with $t-u \le 1$. Here
\begin{eqnarray*}
C^{(k)}_1(a) &=&
\sup_{x\ge a/2,t>0}|\frac{\partial^k g}{\partial x^k}(t,x)|<\infty, \nonumber \\
C^{(k)}_2(q,a)&=&
2a
(\int^{\infty}_0 y|\frac{\partial^k g}{\partial y^k}(1,y)|^q dy)
\sup_{u>0}\frac{g(u,\frac{a}{2})}{u}
<\infty.
\end{eqnarray*}
\end{prop}
{\it Proof.}
We have
$$
\frac{\partial^k g}{\partial x^k}(t,x)
=\frac{\partial^k}{\partial x^k}(t^{-\frac{1}{2}}g(1,t^{-\frac{1}{2}}x))
=t^{-\frac{k+1}{2}}\frac{\partial^k g}{\partial x^k}(1,t^{-\frac{1}{2}}x).
$$
Since $\{B^a_t\}$ and $\{W_t\}$ are independent,
\begin{eqnarray}
&&E[1_{\{\tau^a>u\}}|\frac{\partial^k g}{\partial x^k}(t-u,B^a_u)|^q|\mathcal{G}^W_u]
\nonumber \\
&=&
\int^{\infty}_{0}(g(u,x-a)-g(u,x+a))|\frac{\partial^k g}{\partial x^k}(t-u,x)|^q dx
\nonumber \\
&=&
\int^{\infty}_{0}g(u,x-a)(1-\exp(-\frac{2ax}{u}))
|\frac{\partial^k g}{\partial x^k}(t-u,x)|^q dx
\nonumber \\
&\le&
 \int^{\infty}_{a/2}g(u,x-a)
|\frac{\partial^k g}{\partial x^k}(t-u,x)|^q dx
\nonumber \\
&+&\frac{2a}{u}\int^{a/2}_{0} g(u,x-a)
x|\frac{\partial^k g}{\partial x^k}(t-u,x)|^q dx.\label{equation:g}
\end{eqnarray}
For the first term, we have
$$
\int^{\infty}_{a/2}g(u,x-a)
|\frac{\partial^k g}{\partial x^k}(t-u,x)|^q dx
\le
C^{(k)}_1(a) \int^{\infty}_{a/2} g(u,x-a)dx
\le
C^{(k)}_1(a). 
$$
For the second term, we have
\begin{eqnarray*}
&&\frac{2a}{u}\int^{a/2}_{0}g(u,x-a)x|\frac{\partial^k g}{\partial x^k}(t-u,x)|^q dx
\nonumber\\
&\le&
g(u,\frac{a}{2})\frac{2a}{u}\int^{\infty}_0 x|\frac{\partial^k g}{\partial x^k}(t-u,x)|^q dx 
\nonumber\\
&=&
 g(u,\frac{a}{2}) \frac{2a}{u} \int^{\infty}_0 x|(t-u)^{-\frac{k+1}{2}}\frac{\partial^k g}{\partial x^k}(1,(t-u)^{-\frac{1}{2}}x)|^q dx
\nonumber\\
&= &
g(u,\frac{a}{2})\frac{2a}{u}
(\int^{\infty}_0 y|\frac{\partial^k g}{\partial y^k}(1,y)|^q dy)(t-u)^{\frac{-kq-q+2}{2}}
\nonumber\\
&\le& 
C^{(k)}_2(q,a) (t-u)^{\frac{-kq-q+2}{2}}.
\end{eqnarray*}
Then we have our assertion.
\qed

To represent the conditional expectation under $P$ with respect to $\{\mathcal{G}^W_t\}$and
$\{\mathcal{F}^W_t\}$, we must derive some inequalities to define stochastic integrals.
Propositions \ref{ineq:HkaVer2} and \ref{ineq:Vbara} enable us to evaluate $\bar{H}_a$ and $\bar{U}_a$
in Theorem \ref{th1}. These quantities are defined in Equations (\ref{barH}) and (\ref{barU}), respectively.

\begin{prop}\label{ineq:HkaVer2}
Let $p \in (1, \infty)$ and $q=\frac{p}{p-1}$.
\\
{\rm (1)}
For $k=0,1,2$,  
there are some $C^{(k)}_3(q,a)$ and $C^{(k)}_4(q,a) \in (0, \infty)$ such that
$$
\widetilde{H}^{(k)}_a(t,u;f)
\le
\left(C^{(k)}_3(q,a)
+
C^{(k)}_4(q,a)
(t-u)^{\frac{-kq-q+2}{2q}}\right)
E[|f_u|^p|\mathcal{G}^W_u]^{\frac{1}{p}}
$$
for any $f \in \mathcal{L}^{p}$, $t>u>0$.
Note that $\widetilde{H}^{(k)}_a$ is defined in Equation (\ref{tildeH}).
\\
{\rm (2)} 
Let $k=0,1,2$ and $p>\frac{4}{3-k}$. Then  
there are some $C^{(k)}_{5,1}(q,a)$ and $C^{(k)}_{6,1}(q,a) \in (0,\infty)$
 such that
$$
\int^t_0 \widetilde{H}^{(k)}_a(t,u;f) du
\le
\left(C^{(k)}_{5,1}(q,a) t^{\frac{1}{q}}+C^{(k)}_{6,1}(q,a) t^{\frac{-kq-q+4}{2q}}\right)
\left(\int^t_0 E[|f_u|^p] du\right)^{\frac{1}{p}},
$$
for any $t>0$, $f \in \mathcal{L}^{p}$. 
\\
{\rm (3)} 
Let $k=0,1$, $p>\frac{3}{2-k}$. There are some $C^{(k)}_{5,2}(q,a)$ and $C^{(k)}_{6,2}(q,a) \in (0,\infty)$
 such that
$$
\int^t_0 \widetilde{H}^{(k)}_a(t,u;f)^2 du
\le
\left(C^{(k)}_{5,2}(q,a) t^{\frac{1}{q}} +C^{(k)}_{6,2}(q,a) t^{\frac{-kq-q+3}{q}}\right)
\left(\int^t_0 E[|f_u|^{2p}] du\right)^{\frac{1}{p}},
$$
for any $t>0$, $f \in \mathcal{L}^{2p}$.
 \\
{\rm (4)} 
Let $s \in [0,T]$.
There is some $\hat{C}_1(T,q,a) \in (0,\infty)$ such that
$$
E[\int^s_0 |\bar{H}_a(t;f)| dt] 
\le
 \hat{C}_1(T,q,a)\left(\int^s_0 E[|f_u|^p]du\right)^{\frac{1}{p}},
$$
for any $f \in \mathcal{L}^{p}$, $p > 4$. 
Note that $\bar{H}^{(k)}_a$ is defined in Equation (\ref{barH}).
 \\
 {\rm (5)}
Let $0\le s_0 < s_1$  and 
$\xi$ be a bounded $\mathcal{F}_{s_0}$-measurable random variable.
Then, we have
$$
\int^s_0\hat{H}^{(2)}_a(r;\xi 1_{(s_0,s_1]}(\cdot))dr
=
2\hat{H}^{(0)}_a(s;\xi1_{(s_0,s_1]}(\cdot)).
$$
Note that $\hat{H}^{(k)}_a$ is defined in Equation (\ref{hatH}).
\end{prop}
{\it Proof.}
{\rm (1)} 
By Proposition \ref{ineq:HkaVer2_1}, H\"{o}lder's inequality and a property of convex function, we have
\begin{eqnarray*}
&&|\widetilde{H}^{(k)}_a(t,u;f)|
\nonumber \\
&\le& 
E[1_{\{\tau^a>u\}}|
\frac{\partial^k g}{\partial x^k}(t-u,B^a_u)|^q|\mathcal{G}^W_u]^{\frac{1}{q}}
E[|f_u|^p|\mathcal{G}^W_u]^{\frac{1}{p}}
\nonumber \\
&\le& 
\left(C^{(k)}_1(a)+C^{(k)}_2(q,a) (t-u)^{\frac{-kq-q+2}{2}} \right)^{\frac{1}{q}}.
E[|f_u|^p|\mathcal{G}^W_u]^{\frac{1}{p}}
\nonumber \\
&\le&
\left(C^{(k)}_3(q,a)
+
C^{(k)}_4(q,a)
(t-u)^{\frac{-kq-q+2}{2q}}\right)
E[|f_u|^p|\mathcal{G}^W_u]^{\frac{1}{p}},
\end{eqnarray*}
where
$$
C^{(k)}_3(q,a)=2^{\frac{1}{q}}C^{(k)}_1(a)^{\frac{1}{q}},
\;\;\;
C^{(k)}_4(q,a)=2^{\frac{1}{q}}C^{(k)}_2(q,a)^{\frac{1}{q}}.
$$
Next, we will show assertion (2) and (3).
Let $m = 1,2$.
\begin{eqnarray*}
&&\int^t_0 \widetilde{H}^{(k)}_a(t,u;f)^m du
\nonumber \\
&\le&
\int^t_0
\left(C^{(k)}_3(q,a)+C^{(k)}_4(q,a)(t-u)^{\frac{-kq-q+2}{2q}}\right)^{mq}
E[|f_u|^p|\mathcal{G}^W_u]^{\frac{m}{p}}
du
\nonumber \\
&\le&
\left\{
\int^t_0
\left(C^{(k)}_3(q,a)+C^{(k)}_4(q,a)(t-u)^{\frac{-kq-q+2}{2q}}\right)^{mq}
du
\right\}^{\frac{1}{q}}
\left(\int^t_0
E[|f_u|^p|\mathcal{G}^W_u]^m
du
\right)^{\frac{1}{p}}
\nonumber \\
&\le&
2^m
\left(
\int^t_0
(C^{(k)}_3(q,a)^{mq}
+C^{(k)}_4(q,a)^{mq} (t-u)^{\frac{(-kq-q+2)m}{2}})
du
\right)^{\frac{1}{q}}
\left(\int^t_0
E[|f_u|^p|\mathcal{G}^W_u]^m
du
\right)^{\frac{1}{p}}.
\end{eqnarray*}
If $m=1$ and $p > \frac{4}{3-k}$, or if $m=2$ and $p>\frac{3}{2-k}$, we have $p>\frac{2+2m}{2+m-mk}$ and $\frac{(-kq-q+2)m}{2}>-1$. 
Then we have
$$
\int^t_0
\left(C^{(k)}_3(q,a)^{mq}
+C^{(k)}_4(q,a)^{mq} (t-u)^{\frac{(-kq-q+2)m}{2}}
\right)
du
$$
$$
\le
C^{(k)}_3(q,a)^{mq}t+C^{(k)}_4(q,a)^{mq}\frac{2}{|(-kq-q+2)m+2|}t^{\frac{(-kq-q+2)m+2}{2}}.
$$

Then we have the following for $f \in \mathcal{L}^{mp}$.
$$
\int^t_0 \widetilde{H}^{(k)}_a(t,u;f)^m du
$$
$$
\le
(C^{(k)}_{5,m}(q,a) t^{\frac{1}{q}} +C^{(k)}_{6,m}(q,a) t^{\frac{(-kq-q+2)m+2}{2q}})
(\int^t_0 E[|f_u|^{mp}] du)^{\frac{1}{p}},
$$
where
$$
C^{(k)}_{5,m}(q,a)=2^{\frac{mq+m}{q}}C^{(k)}_1(a)^{\frac{m}{q}},
\;\;\;
C^{(k)}_{6,m}(q,a)=\frac{2^{\frac{mq+m}{q}}}{|(-kq-q+2)m+2|^{\frac{1}{q}}}C^{(k)}_2(q,a)^{\frac{m}{q}}
.
$$
\\
{\rm (4)}
We can see that $H^{(k)}_a(t,f)$, $k=0,1,2,$ are well defined for $f \in \mathcal{L}^{\frac{4}{3-k}+}$ by Assertion (2).
Then $\bar{H}_a(t;f)$ is well defined for $p \in \mathcal{L}^{4+}$.
Since  $p>4$ and ${\frac{-3q+4}{2q}}>0$, 
 Assertion (1) implies
\begin{eqnarray*}
&&E[\int^s_0 |\bar{H}_a(t;f)|dt]
\\
&\le& 
E[e^{\int^s_0 \lambda_a(r) dr}
\int^T_0 
\left(\hat{H}^{(2)}_a(t;f) +2 \lambda_a(t) \hat{H}^{(0)}_a(t;f)\right)dt]
\times 
\left(\int^s_0 E[|f_u|^p]du\right)^{\frac{1}{p}}
\\
&\le& 
\hat{C}_1(T,q,a) \left(\int^s_0 E[|f_u|^p]du\right)^{\frac{1}{p}},
\end{eqnarray*}
where
$$
\hat{C}_1(s,q,a)
$$
$$
=
e^{\int^s_0 \lambda_a(r) dr}
\left\{\int^s_0 
\left(C^{(2)}_3(a)^{\frac{1}{q}}t^{\frac{1}{q}}+C^{(2)}_4(q,a)^{\frac{1}{q}}t^{\frac{-3q+4}{2q}}\right)
+2\lambda_a(t)\left(C^{(0)}_3(a)^{\frac{1}{q}}t^{\frac{1}{q}}+C^{(0)}_4(q,a)^{\frac{1}{q}}t^{\frac{-q+4}{2q}}\right)dt
\right\}.
$$
\\
{\rm (5)}
Since $1_{\{\tau^a>u\}}\int^s_u \frac{\partial g}{\partial r}(r-u,B^a_u)dr=1_{\{\tau^a>u\}}g(s-u,B^a_u)$,
we have the following.
\begin{eqnarray*}
&&\int^s_0\hat{H}^{(2)}_a(r;\xi 1_{(s_0,s_1]}(\cdot))dr
\\
&=&\int^s_0\left(\int^r_0 H^{(2)}_a(r,u;\xi 1_{(s_0,s_1]}(\cdot))du\right)dr
\\
&=&\int^s_0\left(\int^s_u H^{(2)}_a(r,u;\xi 1_{(s_0,s_1]}(\cdot))dr\right)du
\\
&=&\int^s_0\left(\int^s_u E[1_{\{\tau^a>u\}}\xi 1_{(s_0,s_1]}(u)\frac{\partial^2 g}{\partial x^2}(r-u,B^a_u)|\mathcal{G}^W_u]dr\right)du
\\
&=&2\int^s_0\left(\int^s_u E[1_{\{\tau^a>u\}}\xi 1_{(s_0,s_1]}(u)\frac{\partial g}{\partial r}(r-u,B^a_u)|\mathcal{G}^W_u]dr\right)du
\\
&=&2\hat{H}^{(0)}_a\left(s;\xi1_{(s_0,s_1]}(\cdot)\right).
\end{eqnarray*}
Note that the last equation holds by Assertion (2)
\qed
\begin{prop}\label{ineq:Vbara}
Let $T>0$, $p>3$, $q=\frac{p}{p-1}$.
Then  $\bar{U}_a$ is well defined, for any $f \in \mathcal{L}^{6+}$ and
there are 
$\widetilde{C}_1(q,a,T)$, $\widetilde{C}_2(a,T) \in (0, \infty)$
such that
$$
E[\int^T_0 \left(\int^t_0 \bar{U}_a(t,u;f)^2 du\right) dt]
\le 
\widetilde{C}_1(q,a,T)
\left(
\int^T_0 E[|f_u|^{2p}] du
)
dt
\right)^{\frac{1}{p}}
+
\widetilde{C}_2(a,T)
E[\int^T_0 f_u^2 du]
$$
for any $f \in \mathcal{L}^{6+}$.
Note that $\bar{U}$ is given by Equation(\ref{barU}).
\end{prop}
{\it Proof.}
Because $0\le \Phi(t-s,B^a_s) \le 1$, for any $f \in \mathcal{L}^{6+}$ ,we have
\begin{eqnarray*}
&&\int^t_0 E[U_a(t,u;f)^2]du
\\
&\le& \int^t_0 E[E[1_{\{\tau^a>s\}}f_u (2\Phi(t-u,B^a_u)-1)|\mathcal{G}^W_u]^2]du
\\
&\le& \int^t_0 E[ f_u^2 (2 \Phi(t-u,B^a_u)-1)^2]du
\le \int^t_0 f_u^2 du.
\end{eqnarray*}
By the above evaluation and Proposition \ref{ineq:HkaVer2} (2), we have
\begin{eqnarray*}
&&E[\int^T_0 \left(\int^t_0 \bar{U}_a(t,u;f)^2 du \right)dt]
\\
&=&
E[\int^T_0 \left(\int^t_0 e^{2\int^t_0 \lambda_a(r) dr}(H^{(1)}_a(t,u;f)+\lambda_a(t)^2 U_a(t,u;f))^2 du \right) dt]
\\
&\le& 
2e^{2\int^T_0 \lambda_a(r) dr} 
E[\int^T_0 
\int^t_0 
\left(\tilde{H}^{(1)}_a(t,u;f)^2 du
\right)
 dt
+\int^T_0
\lambda_a(t)^2
\left(\int^t_0 |U_a(t,u;f)| du
\right)^2
 dt]
\\
&\le& 
2e^{2\int^T_0 \lambda_a(r) dr} 
\int^T_0 
\left(
(C^{(1)}_{5,2}(q,a) t+C^{(1)}_{6,2}(q,a)t^{-2q+3})^{\frac{1}{q}}
\lambda_a(t)^2
(\int^t_0 E[|f_u|^p]^2 du)^{\frac{1}{p}}
\right)
dt
\\
&+&
2T e^{2\int^T_0 \lambda_a(r) dr} 
\left(\sup_{0 \le t \le T}\lambda_a(t)^2\right)
E[\int^T_0 f_u^2 du].
\end{eqnarray*}
For a part of first term, we have
$$
\int^T_0 (
(C^{(1)}_{5,2}(q,a) t+C^{(1)}_{6,2}(q,a)t^{-2q+3})^{\frac{1}{q}}
(\int^t_0 E[|f_u|^{2p}] du)^{\frac{1}{p}})
dt
$$
$$
\le
(\int^T_0 
(
C^{(1)}_{5,2}(q,a) t+C^{(1)}_{6,2}(q,a)t^{-2q+3}
)
dt)^{\frac{1}{q}}
(\int^T_0 
(
\int^t_0 E[|f_u|^{2p}] du
)
dt
)^{\frac{1}{p}}
.
$$
Note that $U$ is defined in Equation (\ref{U}).
Then we have the assertion where
\begin{eqnarray*}
&&\widetilde{C}_1(q,a,T)\\
&=&
2e^{2\int^T_0 \lambda_a(r) dr} 
\left(\int^T_0 
(
C^{(1)}_{5,2}(q,a) t+C^{(1)}_{6,2}(q,a)t^{-2q+3}
)
dt\right)^{\frac{1}{q}}
\left(\int^T_0 
(
\lambda_a(t)^2
\int^t_0 E[|f_u|^{2p}] du
)
dt
\right)^{\frac{1}{p}}
\end{eqnarray*}
and
$$
\widetilde{C}_2(a,T)
=
2T e^{2\int^T_0 \lambda_a(r) dr}
\left( 
\sup_{0 \le t \le T}\lambda_a(t)^2
\right)
.
$$
\qed

\section{Representation theorem}
We saw that some integrals are well defined under the conditions in Section \ref{sec:2}.
In this section, we prove Theorem \ref{th1}, which is the representation theorem under $\mathcal{F}^W_t$.
 For $x,y \ge 0$ and $t>0$, let
\begin{eqnarray}\label{g_0}
g_0(t,x,y) = g(t,y-x) - g(t,y+x)
= g(t,y-x)(1 - e^{-2xy/t})  
\end{eqnarray}
where
$g(t,x)$ and $\Phi(t,x)$ are the density and distribution, respectively, of the Brownian motion $B_t$. These are given by Equation (\ref{g}).

First, we will present a representation theorem for $E[\int^t_0 \cdot dB_s|\mathcal{F}^W_t]$ which corresponds to Theorem \ref{th1}(1).

\begin{lemma}\label{eq:infty_u}
Let $t>u>0$ and $\xi$ be a bounded $\mathcal{B}_u$-measurable random variable. Then we have
$$
E[\xi|\mathcal{G}^W_{\infty}]=E[\xi|\mathcal{G}^W_{u}]\;\; and \;\; E[\xi(B_t-B_u)|\mathcal{G}^W_{\infty}]=0.
$$
\end{lemma}
{\it Proof.}
Let $h_0$ be a bounded $\mathcal{G}^W_u$ -measurable random variable and 
$h_1$ be a bounded $\sigma\{W(s)-W(u);s \ge u\}$ measurable random variable.Then
\begin{eqnarray*}
&&E[\xi h_0 h_1]
\\
&&=E[\xi h_0 E[h_1|\mathcal{B}_u]]=E[\xi h_0]E[h_1]
\\
&&=E[E[\xi|\mathcal{G}^W_u] h_0]E[h_1]
=E[E[\xi|\mathcal{G}^W_u]h_0h_1]
\end{eqnarray*}
and 
$$
E[\xi (B_t-B_u)h_0h_1]=E[\xi h_0]E[(B_t-B_u)h_1]=0.
$$
So we have our assertion.
\qed
\begin{prop}\label{eq:B0}
Let $0\le s_0 < s_1$, $\xi$ be a bounded $\mathcal{B}_{s_0}$-measurable random variable.
Then, we have the following for 
$t \ge 0$,
$$
E[\xi 1_{\{\tau^a>t\}} (B^a_{s_1}-B^a_{s_0})]
$$
\begin{equation}\label{equation:B}
= 
-
\int^{\infty}_{t} 
(\int^{s_1}_{s_0}
1_{\{ u< r\} }E[\xi  1_{\{\tau^a > u\} }\frac{\partial^2 g}{\partial x^2}(r-u,B^a_u)]du) dr.
\end{equation}
\end{prop}
{\it Proof.}
Let
$$
\varphi (s,x,t)
= 
\int_{0}^{\infty}\int_{0}^{\infty}
(y-x)g_0(s,x,y)g_0(t,y,z) dydz,
\qquad x>0,
\;
s,t>0.
$$
Note that $g_0$ is defined in Equation (\ref{g_0}).
At first, let us think about the case $t>s_1$. Then we have
$$
1_{\{ \tau^a >s_0\} }E[1_{\{ \tau^a >t\}}(B^a_{s_1}-B^a_{s_0})|{\cal B}_{s_0}]
= 1_{\{ \tau^a >s_0\} } \varphi (s_1-s_0,B^a_{s_0},t-s_1).
$$
Then 
$$
E[\xi 1_{\{ \tau^a > t\} } (B^a_{s_1}-B^a_{s_0})]
= E[ \xi  1_{\{ \tau^a >s_0\}} \varphi (s_1-s_0,B^a_{s_0}, t-s_1)].
$$
Note that
$$
|\varphi (s,x,t)|
\le 
\int^{\infty}_{-\infty} \int^{\infty}_{-\infty} |y-x|g(s,x-y)g(t,y-z)dzdy
$$
\begin{equation}\label{eq:bdd}
=\int^{\infty}_{-\infty}|y-x|g(s,x-y)dy = E[|B_s|]
=\sqrt{\frac{2s}{\pi}}.
\end{equation}
Since
$$
d_s\varphi(s_1-s,B^a_s,r)=
(-\frac{\partial}{\partial s}+\frac{1}{2}\frac{\partial^2}{\partial x^2})\varphi(s_1-s,B^a_s,r)ds
+\frac{\partial \varphi}{\partial x}(s_1-s,B^a_s,r)dB^a_s
$$
and
\begin{eqnarray*}
&&(- \frac{\partial}{\partial s}+\frac{1}{2} \frac{\partial^2 }{\partial x^2} )
\varphi (s,x,r)
\\
&=& - \int^{\infty}_{0}\int^{\infty}_{0}
\frac{\partial g_0}{\partial x}(s,x,y)g_0(r,y,z) dydz
= - \int^{\infty}_{0}\frac{\partial g_0}{\partial x}(s+r,x,z) dz
\\
&=& - \int^{\infty}_{0}
( \frac{\partial g}{\partial x}(s+r,x-z)
- \frac{\partial g}{\partial x}(s+r,x+z)) dz
= -2g(s+r,x), \;\;\; x>0, \;\; s,r>0,
\end{eqnarray*}
we have
$$
1_{\{ \tau^a >s_0\}}
( \varphi (s_1-s\wedge \tau^a,B^a_{s\wedge \tau^a }, t-s_1)
- \varphi (s_1-s_0,B^a_{s_0}, t-s_1))
$$
$$
= -2 \int_{s_0}^{s\wedge \tau^a} g(t-u, B^a_u)du 
+ \int_{s_0}^{s\wedge \tau^a} \frac{\partial \varphi }{\partial x}
(s_1-u,B^a_u, t-s_1)dB^a_u,\qquad s\in[s_0,s_1).
$$
As $2\frac{\partial g}{\partial t}=\frac{\partial^2 g}{\partial x^2}$, we have
\begin{eqnarray*}
&&E[\xi 1_{\{ \tau^a > t\} } (B^a_{s_1}-B^a_{s_0})]
\\
&=& E[ \xi  1_{\{ \tau^a >s_0\}} \varphi (s_1-s\wedge \tau^a,B^a_{s\wedge \tau^a}, t-s_1)]
\\
&+&2 E[\xi  1_{\{ \tau^a >s_0\}}
(\int_{s_0}^{s} 1_{\{\tau^a>u \} }g(t-u, B^a_u)du )], \;\;\; s \in [s_0,s_1).
\end{eqnarray*}
Since
$\varphi (s,0,t)=0$ and $\varphi (s,x,t) \to 0,$ $s \downarrow 0,$
we have
$$
\lim_{s \rightarrow s_1} E[\xi  1_{\{ \tau^a >s_0\}} \varphi (s_1-s \wedge \tau^a,B^a_{s \wedge \tau^a}, t-s_1)]
\rightarrow 0
$$
by Equation (\ref{eq:bdd}) and the bounded convergence theorem.
Then we have
\begin{eqnarray*}
&&E[\xi 1_{\{ \tau^a > t\} } (B^a_{s_1}-B^a_{s_0})]
\nonumber
\\
&=&
- 2 \int_{s_0}^{s_1} E[\xi  1_{\{ \tau^a > u\}}g(t-u, B^a_u)]drdu
\nonumber
\\
&=&
-2\int^{s_1}_{s_0} E[\xi \int^{\infty}_t \frac{\partial g}{\partial r}(r-u,B^a_u)dr] du
\nonumber
\\
&=& 
-\int^{\infty}_t (\int_{s_0}^{s_1} 1_{\{u<r\}}E[\xi  1_{\{ \tau^a >u\}}\frac{\partial^2 g}{\partial x^2}(r-u, B^a_u)]du)dr
\end{eqnarray*}
for any $t >s_1$.
 By taking $t \downarrow s_1$, we also have our assertion for $t =s_1$.
\\
Second, let us think of the case $t \in (s_0,s_1]$. 
$$
E[\xi 1_{\{ \tau^a > t\} } (B^a_{s_1}-B^a_{s_0})]
=E[\xi 1_{\{ \tau^a > t\} } E[(B^a_{s_1}-B^a_{s_0})|{\cal B}_{t}]]
$$
$$
=E[\xi 1_{\{ \tau^a > t\} } (B^a_{t}-B^a_{s_0})]
=
-
\int^{\infty}_t 
(
\int^t_{s_0} 1_{\{ u<r \} }E[\xi 1_{\{ \tau^a >u\} }\frac{\partial^2 g}{\partial x^2}(r-u,B^a_u)]
du )
dr
.
$$
Let $p>4$ and $q=\frac{p}{p-1},r>u \ge 0$. Then we have
\begin{eqnarray}
&&E[|1_{\{\tau^a >u\}} \xi \frac{\partial^2 g}{\partial x^2}(r-u,B^a_u)|]
\nonumber\\
&\le&
E[|1_{\{\tau^a >u\}} |\xi|^p]^{\frac{1}{p}}
E[|1_{\{\tau^a >u\}} |\frac{\partial^2 g}{\partial x^2}(r-u,B^a_u)|^q]^{\frac{1}{q}}
\nonumber\\
&\le&
E[|1_{\{\tau^a >u\}} |\xi|^p]^{\frac{1}{p}}
E[(C^{(2)}_1(a)+C^{(2)}_2(2,a)(r-u)^{\frac{-3q+2}{2}}]^{\frac{1}{q}}
\label{bddg2}
\end{eqnarray}
by Proposition \ref{ineq:HkaVer2_1}. 
We have the following by Lemma \ref{eq:infty_u}.
$$
\int^{\infty}_t (\int^{s_1}_{t} 1_{\{ u<r \} }E[1_{\{\tau^a >u\}} \xi \frac{\partial^2 g}{\partial x^2}(r-u,B^a_u)] du)dr
$$
$$
=2\int^{s_1}_{t} ( \int^{\infty}_t 1_{\{ u<r \} }E[1_{\{\tau^a >u\}} \xi \frac{\partial g}{\partial r}(r-u,B^a_u)] dr)du
$$
$$
=-2 \int^{s_1}_{t} ( E[1_{\{\tau^a >u\}} \xi \int^{\infty}_u  \frac{\partial g}{\partial r}(r-u,B^a_u)]dr)du
=0.
$$
Note that since $\frac{-3q+2}{2}>-1$ and by Equation(\ref{bddg2}), we can use Fubini's Theorem in the above equation.
So we have Equation (\ref{equation:B}) for $t \in (s_0,s_1]$.
\\
When $t \in[0,s_0]$,
$$
E[\xi 1_{\{ \tau^a > t\} } (B^a_{s_1}-B^a_{s_0})]
=E[\xi 1_{\{ \tau^a > t\} }E[B^a_{s_1}-B^a_{s_0}|\mathcal{B}_{s_0}]]
=0.
$$
So we see Equation (\ref{equation:B}) is valid for $t \ge 0$.
\qed
\begin{prop}\label{prop_xi_tau}
Let $0\le s_0<s_1$, $t>0$, and $\xi$ be a bounded $\mathcal{F}_{s_0}$-measurable random variable.
Then, we have
$$
E[\xi 1_{\{\tau^a > s_0 \}} 1_{\{\tau^a > t \}}]
=
-
  \int_{s_0 \vee t}^{\infty}  
E[1_{\{ \tau^a > s_0\} } \xi \frac{\partial g}{\partial x} (r-s_0,B^a_{s_0})]dr.
$$
\end{prop}
{\it Proof.}
We assume that $t>s_0$, then we have
$$
E[\xi 1_{\{ \tau^a >t \} }]
= E[\xi 1_{\{ \tau^a > s_0\} }
E[1_{\{ \tau^a > t\} }|{\cal B}_{s_0}]]
= E[\xi 1_{\{ \tau^a > s_0\} } ( \int_{0}^{\infty} g_0(t- s_0,B^a_{s_0},y ) dy)] .
$$
For $x>0$ and $t>0$, we have
\begin{eqnarray*}
&&\int^{\infty}_{0} g_0(t,x,y ) dy 
= - \int^{\infty}_{0} \left(\int^{\infty}_{t} \frac{\partial g_0}{\partial s}(s,x,y ) ds \right) dy 
\\
&=& - \frac{1}{2} \int^{\infty}_{t}
 \left(\int^{\infty}_{0}\frac{\partial^2 g_0}{\partial y^2} (s,x,y) dy\right) ds
= \frac{1}{2} \int^{\infty}_{t}
 \frac{\partial g_0}{\partial y} (s,x,0) ds
= 
-
\int^{\infty}_{t} \frac{\partial g}{\partial x}(s,x) ds.
\end{eqnarray*}
Considering Equation (\ref{equation:g}) in Proposition \ref{ineq:HkaVer2_1}, we have 
$$
E[\xi 1_{\{\tau^a > s_0\}}1_{\{\tau^a > t\}}]
=
-
\int_{s_0 \vee t}^{\infty}  
E[\xi 1_{\{ \tau^a > s_0\} } \frac{\partial g}{\partial x} (r-s_0,B^a_{s_0})]dr.
$$
\qed

\begin{prop}\label{eq:B}
Let $0\le s_0 < s_1$,  $\xi$ be a bounded $\mathcal{F}_{s_0}$-measurable random variable,
 and $v:[0,\infty) \to {\bf R}$ be a bounded Borel measurable function. Then we have the following.
\\
{\rm (1)} 
$$
E[\xi  (B^a_{s_1}-B^a_{s_0})v({\tau^a})]
= 
-
\int_{0}^{\infty} 
v(r) \left(\int_{s_0}^{s_1}
1_{\{ u< r\} }E[\xi  1_{\{\tau^a > u\} }\frac{\partial^2 g}{\partial x^2}(r-u,B^a_u)]du \right) dr.
$$
{\rm (2)} 
$
E[\xi 1_{\{ \tau^a > s_0 \} } v(\tau^a)]
=  
-
\int_{s_0}^{\infty} v(r) 
E[1_{\{ \tau^a > s_0\} } \xi  \frac{\partial g}{\partial x} (r-s_0,B^a_{s_0})]dr.
$\\
{\rm (3)} 
$
E[\xi  (B^a_{s_1}-B^a_{s_0})|\mathcal{F}^W_{\infty}]
=
-
\int^{\infty}_0 \gamma_a^{-1}(r)\hat{H}^{(2)}_a(r;\xi 1_{(s_0,s_1]}(\cdot))dN^a_r.
$\\
{\rm (4)} 
$
E[\xi  (B^a_{s_1}-B^a_{s_0})|\mathcal{F}^W_{t}]
= 
-
\int_{0}^{t}\bar{H}_a(s;\xi 1_{(s_0,s_1]}(\cdot)) \lambda_a(s)^{-1} dM^a_s,\;\;\;t>0.
$
\\
\end{prop}
{\it Proof.}
{\rm (1)}
For $v=1_{[t,\infty)}$, Assertion (1) is valid
by Proposition \ref{eq:B0}.
Let $\mathcal{V}$ be the collection of bounded measurable functions $v$ which satisfy Assertion (1).
Then $\mathcal{V}$ is a vector space.
In addition, if $\{v_n\}_{n \in \mathbf{N}}$ is an increasing sequence of non-negative functions
in $\mathcal{V}$ and if $\lim_{n \rightarrow \infty} v_n$ exists and bounded then $\lim_{n \rightarrow \infty} v_n \in \mathcal{V}$.
Let $\mathcal{A}=\{A \subset \mathbf{R};1_A \in \mathcal{V}\}$ then $(t,\infty)\in A$ for each $t>0$. 
$\mathcal{A}$ is $\pi$-system by the monotone convergence Theorem and 
$\mathcal{A}'=\{(t,\infty);t>0\} \in \mathcal{A}$ is $\pi$-system.
Then we have our assertion by the monotone class theorem.
\\
{\rm (2)}
By the same way with Assertion (1), we see that this assertion 
 is valid for any bounded Borel measurable function $v:(0,\infty )\to {\bf R}$ using Proposition \ref{prop_xi_tau}.
 This completes the proof of Assertion.
\\
{\rm (3)}
Let 
$h_0$ be a bounded ${\cal G}_{s_0}^{W}$-measurable Borel function
and $h_1$ be a bounded $\sigma \{ W_t-W_{s_0};\; t>s_0 \}$-measurable function.
Note that $\mathcal{B}_{s_0} \vee \mathcal{G}^B_{\infty}$ and $\sigma\{W_t-W_{s_0};t>s_0\}$ are independent.
By Lemma \ref{eq:infty_u} and
Proposition \ref{eq:B} (1), we have
\begin{eqnarray*}
&&E[\xi (B^a_{s_1}-B^a_{s_0})h_0h_1]
= 
-
E[h_0\xi (B^a_{s_1}-B^a_{s_0})]E[h_1]
\\
&=&
-
\left(\int^{\infty}_{0} 
 (\int_{s_0}^{s_1}
1_{\{ u< r\} }E[h_0\xi  1_{\{ \tau^a > u\} }
\frac{\partial^2 g}{\partial x^2}(r-u,B^a_u)
]du) dr
\right)
E[h_1]
\\
&=&
-
\left(\int^{\infty}_{0} 
 (\int_{s_0}^{s_1}
1_{\{ u< r\} }E[h_0h_1E[\xi  1_{\{ \tau^a > u\} }
\frac{\partial^2 g}{\partial x^2}(r-u,B^a_u)|\mathcal{G}^W_{\infty}]]du) dr
\right)
\\
&=& 
-
E[h_0 h_1(\int^{\infty}_{0} 
\left (\int^{\infty}_{0}
1_{\{ u< r\} }
E[\xi 1_{(s_0,s_1]}(u) 1_{\{ \tau^a > u\} }
\frac{\partial^2 g}{\partial x^2}(r-u,B^a_u)|\mathcal{G}^W_{\infty}]du]) dr\right)
\\
&=&
-
E[h_0 h_1  \gamma_a(\tau^a)^{-1}
\int^{\infty}_{0}
1_{\{\tau^a > u\} }
E[\xi 1_{(s_0,s_1]}(u) 1_{\{ \tau^a > u\} }
\frac{\partial^2 g}{\partial x^2}(\tau^a-u,B^a_u)|\mathcal{G}^W_{u}]du]
\\
&=&
-
E[h_0 h_1 \gamma_a(\tau^a)^{-1}
\hat{H}^{(2)}_a(\tau^a;\xi 1_{(s_0,s_1]}(\cdot))]
=
-
E[h_0 h_1 
\int^{\infty}_0 \gamma_a(r)^{-1}
\hat{H}^{(2)}_a(r;\xi 1_{(s_0,s_1]}(\cdot))dN^a_r].
\end{eqnarray*}
Then we have the assertion.
\\
{\rm (4)} 
We note that $\hat{H}^{(2)}_a(t;\xi 1_{(s_0,s_1]}(\cdot))=0$ for $t \le s_0$ 
and 
$$
E[\xi  (B^a_{s_1}-B^a_{s_0})|\mathcal{F}^W_{s_0}]
=E[E[\xi  (B^a_{s_1}-B^a_{s_0})|{\cal B}_{s_0}]|\mathcal{F}^W_{s_0}]
=0.
$$
Then we have
$$
E[\xi  (B^a_{s_1}-B^a_{s_0})|\mathcal{F}^W_{t}]
=0
= \int_{0}^{t}\bar{H}_a(s;\xi 1_{(s_0,s_1]}(\cdot))\lambda_a(s)^{-1} dM^a_s, \;\;\;\; t \le s_0.
$$
By Lemma \ref{eq:infty_u}, we have
$$
\int^{\infty}_{0} \hat{H}^{(2)}_a(s;\xi 1_{(s_0,s_1]}(\cdot))ds
=E[\xi (B^a_{s_1}-B^a_{s_0})|\mathcal{G}^W_{\infty}]
= 0
$$
and then
\begin{equation}\label{eq:intHhat2}
E[\int_{t}^{\infty} \hat{H}^{(2)}_a(s;\xi 1_{(s_0,s_1]}(\cdot))ds |{\cal G}_{t}^{W}]
=-\int_{0}^{t} \hat{H}^{(2)}_a(s;\xi 1_{(s_0,s_1]}(\cdot))ds.
\end{equation}
By Assertion (3) and Equation (\ref{eq:intHhat2}), we see that
\begin{eqnarray}
&&E[\xi  (B^a_{s_1}-B^a_{s_0})|\mathcal{F}^W_{t}]
\nonumber
\\
&=&
-
E[\int^{\infty}_0 \gamma_a^{-1} (r) \hat{H}^{(2)}_a(r;\xi 1_{(s_0,s_1]}(\cdot))dN^a_r|\mathcal{F}^W_t]
\nonumber
\\
&=&
-
\int^{t}_0 \gamma_a^{-1} (r) \hat{H}^{(2)}_a(r;\xi 1_{(s_0,s_1]}(\cdot))dN^a_r
-
E[\int^{\infty}_t \gamma_a^{-1} (r) \hat{H}^{(2)}_a(r;\xi 1_{(s_0,s_1]}(\cdot))dN^a_r|\mathcal{F}^W_t]
\nonumber
\\
&=&
-
\int_0^t e^{\int_{0}^{s}\lambda_a(r)dr} \hat{H}^{(2)}_a(s;\xi 1_{(s_0,s_1]}(\cdot))\lambda_a(s)^{-1}dM^a_{s}
\nonumber
\\
&=&
-
\int_0^t e^{\int_{0}^{s}\lambda_a(r)dr} \hat{H}^{(2)}_a(s;\xi 1_{(s_0,s_1]}(\cdot))(1-N^a_s)ds
\nonumber
\\
&+&
e^{\int_{0}^{t}\lambda_a(r)dr}(1 - N^a_t) \int_{0}^{t} \hat{H}^{(2)}_a(s;\xi 1_{(s_0,s_1]}(\cdot))ds.\label{eq:eNH2}
\end{eqnarray}
Note that $\gamma_a$ is defined in Equation (\ref{gamma}).
We also note that 
$e^{\int^{t}_{0} \lambda_a(r)dr}(1-N^a_t) = 1 -  \int^{t}_{0} e^{\int^{s}_{0} \lambda_a(r) dr} dM^a_s$ .
We now see that
\begin{eqnarray*}
&&
e^{\int^{t}_{0}\lambda_a(r)dr}(1-N^a_t)\int^t_0 \hat{H}^{(2)}_a(s;\xi 1_{(s_0,s_1]}(\cdot))ds
\\
&=&
\int^t_0 \hat{H}^{(2)}_a(s;\xi 1_{(s_0,s_1]}(\cdot))ds
-
\left(\int^{t}_{0} e^{\int^{s}_{0} \lambda_a(r) dr} dM^a_s\right)\left(\int^t_0 \hat{H}^{(2)}_a(s;\xi 1_{(s_0,s_1]}(\cdot))ds\right)
\\
&=&
-\int^{t}_{0} e^{\int^{s}_{0} \lambda_a(r) dr} \left(\int^s_0 \hat{H}^{(2)}_a(r;\xi 1_{(s_0,s_1]}(\cdot))dr\right) dM^a_s
\\
&+&
\int^t_0 
\hat{H}^{(2)}_a(s;\xi 1_{(s_0,s_1]}(\cdot))
ds
+
\int^t_0 
\left(\hat{H}^{(2)}_a(s;\xi 1_{(s_0,s_1]}(\cdot))
\int^{s}_{0} e^{\int^{r}_{0} \lambda_a(u) du} dM^a_r
\right)ds
.
\end{eqnarray*}
Then, we have the following for $t\ge s_0$,
$$
E[\xi  (B^a_{s_1}-B^a_{s_0})|\mathcal{F}^W_{t}]
=
-
\int_0^t e^{\int_{0}^{s}\lambda_a(r)dr}\hat{H}^{(2)}_a(s;\xi 1_{(s_0,s_1]}(\cdot))\lambda_a(s)^{-1}dM^a_s
$$
$$
+
\int_{0}^{t} 
(\int_{0}^{s} \hat{H}^{(2)}_a(r;\xi 1_{(s_0,s_1]}(\cdot))dr )
d ( e^{\int_{0}^{s}\lambda_a(r)dr}(1 - N^a_{s-})) 
$$
$$
=
-
\int^t_0 
e^{\int^s_0 \lambda_a(r) dr}
\left(
\hat{H}^{(2)}_a(s;\xi 1_{(s_0,s_1]}(\cdot))+ \lambda_a(s) \int^s_0 \hat{H}^{(2)}_a(r;\xi 1_{(s_0,s_1]}(\cdot))dr)
\right)
 \lambda_a(s)^{-1} dM^a_s.
$$
Finally, we have Assertion by Proposition \ref{ineq:HkaVer2} (4).
\qed

Let $\widetilde{\mathcal{L}}_0$ be the space of progressively measurable processes $\varphi_t$ for which there exist
$\mathcal{B}_{s_k}$- measurable bounded random variables $\xi_{s_k}$  such that
$$
\varphi_t=\sum^{n-1}_{k=0} \xi_{s_k}1_{(s_k,s_{k+1}]}(t),\;\;\; t \ge 0,
$$
where $0\le s_0 < s_1< \cdots <s_n \le T$.
For any $p \ge 1$ and $f \in \mathcal{L}^p$, there exist $f_n \in \widetilde{\mathcal{L}}^0$, $n=1,2,\cdots$, such that
$$
\lim_{n \rightarrow \infty} E[\int^T_0 |f_n(s,\omega)-f(s,\omega)|^p ds ]=0 \;\;\; for\;any\; T>0.
$$
The following gives Theorem \ref{th1}(1).
\begin{col}\label{eq:dB}
Let $T>0$. Then we  have
$$
E[\int^{T}_{0}f_sdB_s|\mathcal{F}^W_{\infty}]
= 
-
\int^{\infty}_{0} 
\gamma_a(s)^{-1}
\left(\int_{0}^{\infty}H^{(2)}_a(s,u;f1_{[0,T]}(\cdot))du\right)dN^a_s
$$
for any $f \in \mathcal{L}^{4+}$ and 
$$
E[\int^{T}_{0}f_sdB_s|\mathcal{F}^W_{t}]
= 
-
\int^{t}_{0} \bar{H}_a(s,f1_{[0,T]}(\cdot))
\lambda_a(s)^{-1}
dM^a_s, \;\;\; t>0
$$
for any $f \in \mathcal{L}^{4+}$. 
\end{col}
{\it Proof.}
Let $s_1>s_0\ge 0$
and ${\widetilde f}$ be a bounded $\mathcal{B}_{s_0}$ -measurable function
and $f_t = {\widetilde f}1_{(s_0,s_1]}(t).$
Then we see that the first and second assertion are valid for $f \in \widetilde{\mathcal{L}}^{(0)}$
 by Proposition \ref{eq:B} (3) and (4), respectively. 
 We can see that 
 $\int_{0}^{\infty}H^{(2)}_a(s,u;f1_{[0,T]}(\cdot))du$ in the first assertion 
is well defined for any $f \in \mathcal{L}^{4+}$
by Proposition \ref{ineq:HkaVer2} (2).
As for the second assertion, 
let us take $\{\widetilde{\xi}_n\} \in \widetilde{\mathcal{L}}_0$ such that
$$
\lim_{n\rightarrow \infty}E[|\widetilde{\xi}_n(r)-f_r|]=0 \; for \; all \; r>0.
$$
Then we have 
$$
E[\int^{T}_{0}\widetilde{\xi}_n(s)dB_s|\mathcal{F}^W_{t}]
= 
-
\int^{t}_{0} \bar{H}_a(s,\widetilde{\xi}_n1_{[0,T]}(\cdot))
\lambda_a(s)^{-1}
dM^a_s, \;\;\; t>0,
$$
by Proposition \ref{eq:B} (4).
Since $\sigma\{W_t;t\ge0\}$ and $\sigma\{N_t;t\ge0\}$ are independent,
we have
\begin{eqnarray*}
&&E[\int^T_0 | (\bar{H}_a(s;\widetilde{\xi}_n)-\bar{H}_a(s; f))| \lambda_a(s)^{-1} dN^a_s]
\\
&=& 
E[\int^T_0 E[| (\bar{H}_a(s;\widetilde{\xi}_n)-\bar{H}_a(s; f))|] \lambda_a(s)^{-1} dN^a_s]
\\  
&=& \int^T_0 E[|(\bar{H}_a(s;\widetilde{\xi}_n- f)|]e^{-\int^s_0 \lambda_a(u) du} ds
\rightarrow 0,\;\;\; as \;\;\;n \rightarrow \infty,\; for \; all \; T>0
\end{eqnarray*}
by Proposition \ref{ineq:HkaVer2} (4)
. 
So $\int^{t}_{0} \bar{H}_a(s,f1_{[0,T]}(\cdot))\lambda_a(s)^{-1} dM_s$ 
 is  well defined, and we have the assertion.
\qed

Second, we will state a representation theorem for $E[\int^t_0 \cdot ds|\mathcal{F}^W_t]$, which corresponds to Theorem \ref{th1}(2).
\begin{prop}\label{eq:U}
Let $s > 0$ and
$f$ be a bounded $\mathcal{F}$-progressively measurable process.
Then, we have the following.
\\
{\rm (1)}
$
E[f_{s} |\mathcal{F}^W_{\infty}]
= E[f_{s} |\mathcal{F}^W_{s}]1_{\{ \tau^a \le s \} }
- 
\int^{\infty}_{s} \gamma_a(r)^{-1}H^{(1)}_a(r,s;f)dN^a_{r}.
$
{\rm (2)} 
$
E[f_{s} |\mathcal{F}^W_{t}]
= 
E[f_{s} |\mathcal{F}^W_{s}]
-
\int^{t}_{s}
\gamma_a(r)^{-1}  \bar{U}_a(r,s;f)
\lambda_a(r)^{-1}
dM^a_r, \;\;\;t>s.
$
\end{prop}
{\it Proof.}
\\
{\rm (1)}
Let 
$s>0$, $h_0$ be a bounded ${\cal G}_{s}^{W}$ -measurable Borel function,
 $h_1$ be a bounded $\sigma \{ W_t-W_{s};\; t>s \}$-measurable function
  and $v:[0,\infty) \rightarrow \mathbf{R}$ be a bounded Borel measurable function.
Then we have 
$$
E[f_{s}v(\tau^a)h_0h_1]
= E[f_{s}1_{\{ \tau^a \le s \} }v(\tau^a)h_0h_1]
+E[h_1] E[f_{s}1_{\{ \tau^a > s \} }v(\tau^a)h_0] 
$$
and
\begin{eqnarray*}
&&E[f_{s}1_{\{ \tau^a \le s \} }v(\tau^a)h_0h_1]
= E[h_1]E[f_{s} 1_{\{ \tau^a \le s \} } v(\tau^a)h_0] 
\\
&&= E[h_1]E[E[f_{s}|\mathcal{F}^W_{s}]1_{\{ \tau^a \le s \} }v(\tau^a)h_0]
= E[E[f_{s}|\mathcal{F}^W_{s}]1_{\{ \tau^a \le s \} }v(\tau^a)h_0h_1].
\end{eqnarray*}
Since $\sigma\{W_t;t\ge0\}$ and $\sigma\{N_t;t\ge0\}$ are independent, we have the following by Proposition \ref{eq:B} (2),
\begin{eqnarray*}
&&E[h_1]E[f_{s}1_{\{ \tau^a > s \} }v(\tau^a)h_0] 
\\
&=&
E[h_1] \int^{\infty}_{s} v(r) E[1_{\{ \tau^a > s\} }  f_{s} 
\frac{\partial g}{\partial x} (r-s,B^a_{s})h_0]dr 
\\
&=&
-
E[h_1] 
\int_{s}^{\infty} v(r) \gamma_a(r) E[\gamma_a(r)^{-1} H^{(1)}_a(r,s;f) h_0 ]dr  
\\
&=&
-
\int_{s}^{\infty} v(r) \gamma_a(r) E[\gamma_a(r)^{-1}H^{(1)}_a(r,s;f) h_0h_1 ]dr  
\\
&=&
-
E[ E[v(r) \gamma_a(r)^{-1}H^{(1)}_a(r,s;f) h_0h_1 ]|_{r=\tau^a }1_{\{\tau^a > s\}}]  
\\
&=&
-
E[\gamma_a(\tau^a)^{-1}H^{(1)}_a(\tau^a,s;f)1_{ \{ \tau^a >s\}}v(\tau^a)h_0h_1] 
\\
&=&
-
E[(\int^{\infty}_{s}\gamma_a(r)^{-1}H^{(1)}_a(r,s;f)dN^a_r) v(\tau^a)h_0h_1].
\end{eqnarray*}
So we have
$$
E[f_{s}1_{\{ \tau^a > s \} }|\mathcal{F}^W_{\infty}]
=
-
\int^{\infty}_{s} \gamma_a(r)^{-1}H^{(1)}_a(r,s;f)1_{ \{ \tau^a >s\} }dN^a_r.
$$
Thus we have Assertion.\\
{\rm (2)}
Note that
$$
\frac{\partial}{\partial t}(2 \Phi(t-s,x)-1)
$$
$$
=2\int^x_{-\infty} \frac{\partial g}{\partial t}(t-s,y)dy = \int^x_{-\infty} \frac{\partial^2 g}{\partial y^2}(t-s,y)dy
=\frac{\partial g}{\partial x}(t-s,x)
$$
and that
$$
2 \Phi(t-s,x)-1
=
-\int^{\infty}_t \frac{\partial}{\partial r}(2 \Phi(r-s,x)-1)dr
= -\int^{\infty}_t \frac{\partial g}{\partial x}(r-s,x)dr.
$$
Here we note that
$$
\lim_{t \rightarrow \infty} \Phi(t-s,x) = \frac{1}{2},\;\;\;x>0
$$
and
$$
\lim_{t \downarrow s} \Phi(t-s,x) = 1, \;\;\;x>0.
$$
Let
$$
L_t = 1- \exp(\int^t_0 \lambda_a(s) ds)(1-N^a_t).
$$
Then we have
$$
dL_t = \exp(\int^t_0 \lambda_a(s) ds)(dN^a_t-\lambda_a(t) (1-N^a_t) dt)
=\exp(\int^t_0 \lambda_a(s) ds)dM^a_t.
$$
We note that
$$
dN^a_t = \exp(-\int^t_0 \lambda_a(s) ds)dL_t + \lambda_a(t) (1-N^a_t) dt
$$
and
\begin{eqnarray*}
\gamma_a(t)^{-1}dN^a_t 
&=& \lambda_a(t)^{-1}dL_t + \exp(\int^t_0 \lambda_a(s) ds) (1-N^a_t) dt
\\
&=& \lambda_a(t)^{-1}dL_t -L_t dt + \exp(\int^t_0 \lambda_a(s) ds)dt.
\end{eqnarray*}
Then we have
\begin{eqnarray*}
U_a(t,s,f)
&=&E[1_{\{\tau^a>s\}} f_s(2\Phi(t-s,B^a_s)-1)|\mathcal{G}^W_s]
\\
&=&-E[\int^{\infty}_t  1_{\{\tau^a>s\}} f_s\frac{\partial g}{\partial x}(r-s, B^a_s)|\mathcal{G}^W_s]dr
\\
&=&
-\int^{\infty}_t H^{(1)}_a(r,s;f)dr.
\end{eqnarray*}
It is obvious that
$$
\int^{\infty}_s H^{(1)}_a(r,s;f)\gamma_a(r)^{-1}dN^a_r
$$
$$
=\int^{\infty}_s H^{(1)}_a(r,s;f)\lambda_a(r)^{-1}dL_r
-\int^{\infty}_s H^{(1)}_a(r,s;f)L_r dr
+\int^{\infty}_s H^{(1)}_a(r,s;f)e^{\int^r_0 \lambda_a(u)}dr.
$$
Here we note that the third term at the last equation is $\mathcal{F}_s$-measurable.
And the second term of the above can be described in the following.
\begin{eqnarray*}
&&-\int^{\infty}_s H^{(1)}_a(r,s;f)L_r dr
\\
&=&
-\int^{\infty}_s H^{(1)}_a(r,s;f)(\int^r_sdL_u+L_s) dr
\\
&=&
-\int^{\infty}_s\left(\int^{\infty}_0 H^{(1)}_a(r,s;f)dr\right)dL_u
-\int^{\infty}_s H^{(1)}_a(r,s;f)L_s dr
\\
&=&
\int^{\infty}_s U_a(r,s;f)dL_r+L_s U_a(s+,s;f).
\end{eqnarray*}
Here we note that the second term at the last equation is $\mathcal{F}_s$-measurable.
Then we have
$$
E[f_s | \mathcal{F}^W_{\infty}]
$$
$$
=\left(
E[f_s|\mathcal{F}^W_s]1_{\{\tau^a \le s\}} 
-
L_s U_a(s+,s;f)
-
\int^{\infty}_s H_a^{(1)}(r,s;f)e^{\int^r_0 \lambda_a(s) ds} dr
\right)
$$
$$
-\int^{\infty}_s
\left(
H^{(1)}_a(r,s;f)\lambda_a(r)^{-1}+U_a(r,s;f)
\right)
dL_r.
$$
The first three terms are $\mathcal{F}^W_s$-measurable and the summation should be equal to $E[f_s|\mathcal{F}^W_s]$. 
The last term is equal to
$$
\int^{\infty}_s \left\{
\exp\left(\int^r_0 \lambda_a(u)du\right)
\left(H^{(1)}_a(r,s;f) + \lambda_a(r) U_a(r,s;f)\right)\lambda_a(r)^{-1} 
\right\}
dM^a_r
$$
$$
=\int^{\infty}_s \bar{U}_a(r,s;f) \lambda_a(r)^{-1} dM^a_r
.
$$
Then we have our assertion.
\qed
\\
The following gives Theorem \ref{th1}(2).
\begin{prop}\label{eq:fH}
Let $T,t>0$ and $f \in \mathcal{L}^{2+}$. Then, we have
$$
E[\int^t_0 f_sds|\mathcal{F}^W_t]
=\int^t_0 E[f_s|\mathcal{F}^W_t] ds
$$
$$
=\int^t_0 E[f_s1_{[0,T](s)}|\mathcal{F}^W_s]ds
-
\int^t_0 \left\{ \left(\int^r_0 \bar{U}_a(s,r;f1_{[0,T]}) ds\right)\lambda_a(s)^{-1} \right\}dM^a_r
.
$$
\end{prop}
{\it Proof.}
Remember that $U_a(t,s;f) =E[1_{\{\tau^a>s\}}f_s (2\Phi(t-s,B^a_s)-1)|\mathcal{G}^W_s]$.
We can see that 
$\int^t_0 (\int^r_0 \bar{U}_a(s,r;f1_{[0,T]}) ds)\lambda_a(s)^{-1}dM^a_r$
is well defined for any $f \in \mathcal{L}^{2+}$
by Proposition \ref{ineq:HkaVer2} (2)
. 
Then we have the assertion by Proposition \ref{ineq:Vbara} and \ref{eq:U}.
\qed

Third, we prove Theorem \ref{th1} (3) as follows.
\begin{prop}\label{eq:W1}
Let $s_1 > s_0\ge  0$,  
and $\xi$ be
 a bounded $\mathcal{F}$ -measurable process.
Then we have 
\\
$$
E[\int^{\infty}_0 \xi 1_{(s_0,s_1]}(r) dW_r|\mathcal{F}^W_{\infty}]
$$
$$
= \int^{\infty}_{0} E[\xi 1_{(s_0,s_1]}(r) |\mathcal{F}^W_r]dW_r
- 
\int^{\infty}_{0} 
\left\{
\left(\int^{r}_{0} \bar{U}_a(r,u;\xi 1_{(s_0,s_1]}(\cdot))dW_u \right)\lambda_a(r)^{-1}
\right\}
dM^a_r.
$$
\\
In particular, for any $T>0$ and  $f \in \mathcal{L}^{6+}$,
$$
E[\int^{\infty}_0 f_r 1_{[0,T]}(r) dW_r|\mathcal{F}^W_{t}]
$$
$$
= \int^{t}_{0} E[ f_r 1_{[0,T]}(r)  |\mathcal{F}^W_r]dW_r
- \int^{t}_{0} 
\left\{
\left(\int^{r}_{0} \bar{U}_a(r,u;f 1_{[0,T]}(\cdot))dW_u \right)\lambda_a(r)^{-1}
\right\}
dM^a_r.
$$
\end{prop}
{\it Proof.}
\\
Note that 
$$
E[\int^{\infty}_0 \xi 1_{(s_0,s_1]}(r) dW_r|\mathcal{F}^W_{\infty}]
= E[\xi |\mathcal{F}^W_{\infty}] (W_{s_1}-W_{s_0}).
$$
By Proposition \ref{eq:U}, we have
$$
E[\xi |\mathcal{F}^W_{\infty}] 
= E[\xi |\mathcal{F}^W_{s_0}] 
-
\int^{\infty}_{s_0} \bar{U}_a(r,s_0;\xi1_{[s_0,\infty)}(\cdot))\lambda_a(r)^{-1}dM^a_r 
$$
and then
\begin{eqnarray*}
&&E[\int^{\infty}_0 \xi 1_{(s_0,s_1]}(r) dW_r|\mathcal{F}^W_{\infty}]
\\
&=& \int_{s_0}^{s_1}E[\xi |\mathcal{F}^W_{r}]dW_r
-
\int^{\infty}_{s_0}  
\left\{
\bar{U}_a(r,s_0;\xi 1_{[s_0,s_1)}(\cdot))(W_{r\wedge s_1}-W_{s_0})\lambda_a(r)^{-1}
\right\}
dM^a_r
\\
&=&
\int^{\infty}_{0}E[\xi 1_{(s_0,s_1]}(r) |\mathcal{F}^W_{r}]dW_r
-
\int^{\infty}_{s_0}
(\int^{r \wedge s_1}_{0} 
\left\{
\bar{U}_a(r,s_0;\xi 1_{(s_0,s_1]}(\cdot))dW_u)
\lambda_a(r)^{-1}
\right\}
dM^a_r.
\end{eqnarray*}
Here we note that 
$$
1_{\{ \tau^a > s_0\} } 
( \frac{\partial g}{\partial x} (t-s\wedge \tau^a ,B^a_{s\wedge \tau^a})
- \frac{\partial g}{\partial x} (t-s_0,B^a_{s_0}))
$$
$$
=1_{\{ \tau^a > s_0\} } 
\int_{s_0}^{s\wedge \tau^a }
\frac{\partial^2 g}{\partial x^2} (t-r,B^a_r)dB^a_r,
\qquad s \in (s_0,t).
$$
Then we have
\begin{eqnarray*}
H^{(1)}_a(t,s_0;\xi 1_{[s_0,s_1)}(\cdot))
&=& E[1_{\{ \tau^a > s_0\} }\xi
1_{[s_0,s_1)}(\cdot)
( \frac{\partial g}{\partial x} (t-s\wedge \tau^a ,B^a_{s\wedge \tau^a}))
|{\cal G}_{s_0}^{W}]
\\
&=&
H^{(1)}_a(t,s;\xi 1_{[s_0,s_1)}(\cdot)),\;\;\; s\in (s_0, t \wedge s_1).
\end{eqnarray*}
Also we have
$$
1_{\{ \tau^a > s_0\} } 
(\Phi(t-s\wedge \tau^a ,B^a_{s\wedge \tau^a})
- \Phi(t-s_0,B^a_{s_0}))
$$
$$
=1_{\{ \tau^a > s_0\} } 
\int_{s_0}^{s\wedge \tau^a }
g (t-r,B^a_r)dB^a_r,
\;\;\; s \in (s_0,t \wedge s_1).
$$
Thus we have
$$
U_a(t,s_0;\xi 1_{[s_0,t \wedge s_1)}(\cdot))
=U_a(t,s;\xi 1_{[s_0,t \wedge s_1)}(\cdot)),\;\;\; s\in (s_0,t \wedge s_1).
$$
Since
$$
\bar{U}_a(r,u;\xi 1_{(s_0,s_1]}(\cdot))=0,\;\;\; r\in [0,s_0],
$$
 we can see that
$\int^{\infty}_{0} (\int^{r}_{0} \bar{U}_a(r,u;\xi 1_{(s_0,s_1]}(\cdot))dW_u)\lambda_a(r)^{-1}dM^a_r$
is well defined.
Then we have the first assertion.
For $\widetilde{\xi} \in \widetilde{\mathcal{L}}_0$, we have the following by the first assertion,
$$
E[\int^{T}_{0} \widetilde{\xi} _r dW_r|\mathcal{F}^W_{\infty}]
= \int^{T}_{0} E[\widetilde{\xi}_r|\mathcal{F}^W_r]dW_r 
-
\int^{T}_{0} \left\{
 \left(\int_{0}^{r} \bar{U}_a(r,u;\widetilde{\xi} )dW_u\right) \lambda_a(r)^{-1} 
 \right\}
 dM^a_r.
$$
Let us take $\{\widetilde{\xi}_n\} \in \widetilde{\mathcal{L}}_0$ such that
$$
\lim_{n\rightarrow \infty}E[\int^T_0 |\widetilde{\xi}_n(r)-f_r| dr]=0 \; for \; all \; T>0.
$$
Since $\sigma\{W_t;t\ge0\}$ and $\sigma\{N_t;t\ge0\}$ are independent,
we have
\begin{eqnarray*}
&&E[\int^T_0 |\int^r_0 (\bar{U}_a(r,u;\widetilde{\xi}_n)-\bar{U}_a(r,u; f)) dW_u| \lambda_a(r)^{-1} dN^a_r]\\
&=& E[\int^T_0 E[|\int^r_0 (\bar{U}_a(r,u;\widetilde{\xi}_n)-\bar{U}_a(r,u; f)) dW_u|] \lambda_a(r)^{-1} dN^a_r]\\
&=& \int^T_0 E[|\int^r_0 (\bar{U}_a(r,u;\widetilde{\xi}_n - f) dW_u|] q_{a}(r) dr\\
&\le& \int^T_0 E[ \int^r_0 (\bar{U}_a(r,u;\widetilde{\xi}_n- f)^2 du]^{1/2} dr\\
&\rightarrow& 0,\;\;\; as \;\;\;n \rightarrow \infty,\; for \; all \; T>0
\end{eqnarray*}
by Proposition \ref{ineq:Vbara} for $f \in \mathcal{L}^{6+}$. 
So we have Assertion.
\qed

Fourth, we show Theorem \ref{th1}(4) as follows.
\begin{prop}\label{eq:Bhat}
Let $T,t>0$ and $\hat{f}^j \in \mathcal{L}^{2+}, j=1,\cdots,d$. Then we have
$$
E[ \int^t_0 \hat{f}_s^j d\hat{B}_s^j|\mathcal{F}^W_t]=0\; , j=1,\cdots,d.
$$
\end{prop}
{\it Proof.}
Because $B$, $\hat{B}$ and $W$ are independent and $\mathcal{F}^W_t \subset \sigma \{B_s, W_s; s \in [0,\infty)\}$,
$$
E[\sum^{d}_{j=1} \int^t_0 \hat{f}_s^j d\hat{B}_s^j |\mathcal{F}^W_t]
=E[\sum^{d}_{j=1} \int^t_0 \hat{f}_s^j d\hat{B}_s^j]
=0.
$$
\qed
\\

Finally, we state Nakagawa's \cite{Nakagawa} representation theorem using a different expression.
\begin{prop}\label{prop:H2}
Let $f \in \mathcal{L}^{4+}$. Then we have
$$
\hat{H}^{(2)}_a(t;f)
=\lim_{u\uparrow t}E[(\int^u_0 1_{\{\tau^a>s\}}f_s dB^a_s)\frac{\partial g}{\partial x}(t-u,B^a_u)|\mathcal{G}^W_t].
$$
The right-hand side of the above corresponds to the representation theorem given by Nakagawa \cite{Nakagawa} .
\end{prop}
{\it Proof.}
Note that
$$
(-\frac{\partial}{\partial u}+\frac{1}{2}\frac{\partial^2}{\partial x^2})\frac{\partial g}{\partial x}(t-u,x)=0
$$
and so
$$
d_u \frac{\partial g}{\partial x}(t-u,B^a_u)
=\frac{\partial^2 g}{\partial x^2}(t-u,B^a_u)dB^a_u, \;\;\; u<t.
$$
By Ito formula, we have
$$
d_u ((\int^u_0 1_{\{\tau^a>s\}}f_s dB^a_s)\frac{\partial g}{\partial x}(t-u,B^a_u))
$$
$$
= 1_{\{\tau^a>u\}}f_u \frac{\partial g}{\partial x}(t-u,B^a_u)dB^a_u
+(\int^u_0 1_{\{\tau^a>s\}}f_s dB^a_s)\frac{\partial^2 g}{\partial x^2}(t-u,B^a_u)dB^a_u
$$
$$
+ 1_{\{\tau^a>u\}}f_u \frac{\partial^2 g}{\partial x^2}(t-u,B^a_u)du, \;\;\; u<t.
$$
And then
\begin{eqnarray*}
&&E[(\int^u_0 1_{\{\tau^a>s\}} f_s dB^a_s)\frac{\partial g}{\partial x}(t-u,B^a_u)|\mathcal{G}^W_t]\\
&=&\int^u_0 E[ 1_{\{\tau^a>s\}} f_s \frac{\partial^2 g}{\partial x^2}(t-s,B^a_s)|\mathcal{G}^W_t]ds\\
&=&\int^u_0 E[ 1_{\{\tau^a>s\}} f_s \frac{\partial^2 g}{\partial x^2}(t-s,B^a_s)|\mathcal{G}^W_s]ds\\
&=&\int^u_0 H^{(2)}(t,s;f)ds, \;\;\; u<t.
\end{eqnarray*}
So we have
$$
\hat{H}^{(2)}_a(t;f)
=\lim_{u\uparrow t} E[(\int^{u}_0 1_{\{\tau^a>s\}}f_s dB^a_s)\frac{\partial g}{\partial x}(t-u,B^a_u)|\mathcal{G}^W_t].
$$
\qed
\section{Equivalent probability measures}\label{section:measure}
We now state a representation theorem for a filtering model with first-passage-type stopping time.
Note that $I$, $F$ are defined in Equations (\ref{I}) and (\ref{F}). 
Operators $\widetilde{D}_0$,$\widetilde{D}_1$,$\widetilde{D}_2$ and $\widetilde{L}$ are defined in Equations (\ref{operators}).
As we defined in Equation (\ref{rho}), let
$$
\rho_t=1+\int^t_{0+} \rho_{s-}(b_0(s,X_s,Z_s)d\widetilde{B}_s + \beta (s,X_{s \wedge \tau},Y_s) d\widetilde{W}_s).
$$
Let $F \in \Sigma$ be given by Equation (\ref{F}) in the Introduction. Then we have
$$
F_t = F_0 + \int^t_0 \left(f_1(s)-\beta(s,X_s,Y_s) f_2(s) -b_0(s,X_s,Z_s) f_3(s)\right) ds
$$
$$
+\int^t_0 f_2(s)d\widetilde{W}_s+ \int^t_0 f_3(s) d\widetilde{B}_s + \int^t_0 f_4(s) d\hat{B}_s
$$
and so
\begin{eqnarray*}
\rho_{t} F_{t\wedge \tau} &=&\rho_0 F_0
+\int^{t}_0 F_{s \wedge \tau} d\rho_s
+\int^{t \wedge \tau}_0 \rho_{s-} dF_s+[\rho,F]_{t\wedge \tau}\\
&=&F_0
+\int^t_{0+} \rho_{s-}(\widetilde{D_1} F)_s d\widetilde{B}_s
+\int^t_{0+} \rho_{s-}(\widetilde{D_2} F)_s d\hat{B}_s\\
&+&\int^t_{0+} \rho_{s-}(\widetilde{D_0} F)_s d\widetilde{W}_s
+\int^t_{0+} \rho_{s-}(\widetilde{L} F)_s ds.
\end{eqnarray*}
Let
\begin{equation}
V(t,s;f) =\widetilde{E}[\rho_{s-} 1_{\{\tau>s\}} f_s (2\Phi(t-s,X_s)-1) |\mathcal{G}^Y_s], \label{V}
\end{equation}

\begin{equation}
\bar{V}(t,s;f) = e^{\int_{0}^{t} \lambda_{x_0}(r)dr} ( I^{(1)}(t,s;f) +\lambda_{x_0}(t) V(t,s;f)), \label{barV}
\end{equation}

\begin{equation}\label{barI}
\bar{I}(t,s;f) = e^{\int^t_0 \lambda_{x_0}(r) dr}\left(\int^s_0 I^{(2)}(t,u;f)du+2\lambda_{x_0}(t)\int^s_0 I^{(0)}(t,u;f) du\right),
\end{equation}

\begin{equation}\label{hatV}
\hat{V}(r,s;F) = 
\widetilde{\rho}_{r-}^{-1} 
e^{\int^r_0 \lambda_{x_0}(u) du}
(
\hat{V}_1(r,s;F)
+
\lambda_{x_0}(r) \hat{V}_2(r,s;F)
),\;\;\; s\le r
\end{equation}
where 

\begin{equation}
\hat{V}_2(r,s;F)  
=\int^s_0 V(r,u;\widetilde{D_0} F)d\widetilde{W}_u \nonumber\\
+\int^s_0 \left(V(r,u;\widetilde{L} F) +2 I^{(0)}(r,u;\widetilde{D_1} F)\right)du \label{V2}
\end{equation}
for $f \in \mathcal{L}^{6+}$ and $F \in \Sigma$.
Then we have the following by Theorem \ref{th1}.
\begin{equation}\label{eq:rhoF}
\widetilde{E}[\rho_{t} F_{t \wedge \tau}|\mathcal{F}_t]=F_0
\end{equation}
$$
-
\int^{t \wedge \tau}_0
\left\{
\left(\bar{I}(r,r;\widetilde{D_1} F)
+(\int^r_0 \bar{V}(r,u;\widetilde{D}_0 F)d\widetilde{W}_u )
+
 (\int^r_0 \bar{V}(r,u;\widetilde{L} F)du)
\right)
\lambda_{x_0}(r)^{-1}
\right\}
d\widetilde{M}_r 
$$
$$
+
\int^{t \wedge \tau}_0 \widetilde{E}[\rho_{r-}(\widetilde{L} F)_r|\mathcal{F}_r] dr
+\int^{t \wedge \tau}_0 \widetilde{E}[\rho_{r-}
(\widetilde{D_1} F)_r|\mathcal{F}_r]d\widetilde{W}_r
$$
$$
=
F_0
-
\int^t_0 \widetilde{\rho}_{r-}\hat{V}(r,r; F)
\lambda_{x_0}(r)^{-1}
d\widetilde{M}_r
$$
$$
+
\int^t_0  \widetilde{E}[\rho_{r-} (\widetilde{L} F)_r|\mathcal{F}_r] dr
+\int^t_0  \widetilde{E}[\rho_{r-} (\widetilde{D_1} F)_r|\mathcal{F}_r]d\widetilde{W}_r.\nonumber
$$
Here we note that
\begin{eqnarray*}
&&\bar{I}(r;\rho(\widetilde{D_1} F))
+\int^r_0 \bar{V}(r,u;\widetilde{D_0} F)d\widetilde{W}_u 
+\int^r_0 \bar{V}(r,u;\widetilde{L} F)du
\nonumber \\
&=&
e^{\int^r_0 \lambda_{x_0}(u) du}
\left\{
\left(\int^r_0 I^{(1)}(r,u;\widetilde{D_0} F)d\widetilde{W}_u
+\int^r_0 (I^{(2)}(r,u;\widetilde{D_1} F)+I^{(1)}(r,u;\widetilde{L} F))du\right)
\right.
\nonumber \\
&+&
\left.
\lambda_{x_0}(r) 
\left(
\int^r_0 V(r,u;\widetilde{D_0} F)d\widetilde{W}_u
+\int^r_0 
(
V(r,u;\widetilde{L} F) 
+
2\int^r_0 I^{(0)}(r,u;\widetilde{D_1} F)
)du\right)
\right\}
\nonumber \\
&=&\widetilde{\rho}_{r-} \hat{V}(r,r; F).
\end{eqnarray*}
We will show that $\hat{V}(r,r;F)=\hat{V}(r;F)$ and that these can be written without using stochastic integrals by Propositions \ref{eq:V1} and \ref{eq:V2}.
\begin{prop}\label{eq:V1}
Let $T>0$ and $F \in \Sigma$.
Then we have
$$
\hat{V}_1(r,s;F)
=
-\frac{\partial g}{\partial x}(r,x_0)F_0
+I^{(1)}(r,s;F)
,\;\;\; 0<s<r\le T
$$
and we can see that the right-hand side of the above equation can be defined even at $r=s$ by $r\downarrow s$.
Note that $\hat{V}_1$ is defined in Equation (\ref{VandV1}).
\end{prop}
{\it Proof.}
Because $\frac{\partial g}{\partial r}(r,x)-\frac{1}{2}\frac{\partial^2 g}{\partial x^2}(r,x)=0$, we have
$$
\frac{\partial g}{\partial x}(r-s,X_s)= \int^s_0 \frac{\partial^2 g}{\partial x^2}(r-u,X_u)d\widetilde{B}_u.
$$
So we have
\begin{eqnarray*}
&&d( \frac{\partial g}{\partial x}(r-u,X_u) \rho_u F_u)\\
&=& \rho_u((\frac{\partial^2 g}{\partial x^2}(r-u,X_u)F_s+\frac{\partial g}{\partial x}(r-u,X_u)(\widetilde{D_1} F)_u)d\widetilde{B}_u
\\
&+&\frac{\partial g}{\partial x}(r-u,X_u)\rho_u(\widetilde{D_2} F)_u d\hat{B}_u
+\frac{\partial g}{\partial x}(r-u,X_u)\rho_u(\widetilde{D_0} F)_u d\widetilde{W}_u
\\
&+&(\frac{\partial g}{\partial x}(r-u,X_u)(\widetilde{L} F)_u
+\frac{\partial^2 g}{\partial x^2}(r-u,X_u)\rho_u(\widetilde{D_1} F)_u)du).
\end{eqnarray*}
Since $\mathcal{G}^Y_{s}$, $\sigma\{\tilde{B}_u,\hat{B}_u; u \le s\}$
 and  $\sigma\{\widetilde{M}_u; u \le s\}$
 are independent, we have
the following for $r>s$.
\begin{eqnarray*}
&&\widetilde{E}[\frac{\partial g}{\partial x}(r-s,X_s) 1_{\{\tau>s\}} \rho_s F_s|\mathcal{G}^Y_s]
\\
&=&\frac{\partial g}{\partial x}(r,x_0) F_0
+\int^s_0 I^{(1)}(r,u;\widetilde{D_0} F)d\widetilde{W}_u
+\int^s_0 (I^{(2)}(r,u;\widetilde{D_1} F)+I^{(1)}(r,u;\widetilde{L} F))du
\\
&=&\frac{\partial g}{\partial x}(r,x_0)F_0+\hat{V}_1(r,s;F)
.
\end{eqnarray*}
Then we have our assertion.
\qed
\begin{prop}\label{eq:V2}
Let $T>0$ and $F \in \Sigma$. Then we have
$$
\hat{V}_2(r,s;F)
=
-(2\Phi(r,x_0)-1)F_0
+
V(r,s;F)
,\;\;\; 0<s<r\le T
$$
and
$$
\lim_{s \uparrow r} \hat{V}_2(r,s;F)
=
-(2\Phi(r,x_0)-1)F_0
+
\widetilde{E}[1_{\{\tau>r\}}\rho_{r} F_r|\mathcal{G}^Y_r]
.
$$
In particular, $\hat{V}(r,r;F)=\hat{V}(r;F)$.
Note that $\hat{V}(r;F)$ and $\hat{V}(r,s;F)$ are defined in Equation (\ref{VandV1}) and (\ref{hatV}), respectively.
\end{prop}
{\it Proof.}
Because $\frac{\partial \Phi}{\partial r}(r,x)-\frac{1}{2}\frac{\partial^2 \Phi}{\partial x^2}(r,x)=0$
 and $\frac{\partial \Phi}{\partial x}(r,x) = g(r,x)$,
$$
\Phi(r-s,X_s)= \int^s_0 g(r-u,X_u)d\widetilde{B}_u.
$$
So we have
$$
(2\Phi(r-s,X_s)-1)\rho_s F_s 
$$
$$
=
(2\Phi(r,x_0)-1) F_0
+2\int^s_0 \rho_{u-} F_u g(r-u,X_u) d\widetilde{B}_u
$$
$$
+
\int^s_0
\{
(2\Phi(r-u,X_u)-1)
(
\rho_{u-} (\widetilde{D_1} F)_u d\widetilde{B}_u
+\rho_{u-}(\widetilde{D_2} F)_ud\hat{B}_u
+\rho_{u-}(\widetilde{D_0} F)_u d\widetilde{W}_u
+\rho_{u-}(\widetilde{L} F)_u du
)
\}
$$
$$
+2\int^s_0\rho_{u-}g(r-u,X_u)(\widetilde{D_1} F)_udu
.
$$
And then we have the following by Lemma \ref{eq:infty_u}, which gives the first assertion.
\begin{eqnarray*}
&&V(r,s;F)\\
&=&\widetilde{E}[(2\Phi(r,x_0)-1)F_0
+\int^s_0 1_{\{\tau > u\}}(2\Phi(r-u,X_u)-1)\rho_{u-}(\widetilde{D_0} F)_u d\widetilde{W}_u
\\
&+&\int^s_0 1_{\{\tau > u\}}((2\Phi(r-u,X_u)-1)\rho_{u-}(\widetilde{L} F)_u
+2g(r-u,X_u)\rho_{u-}(\widetilde{D_1} F)_u)du|\mathcal{G}^Y_r]
\\
&=&
(2\Phi(r,x_0)-1) F_0
+\int^s_0 V(r,u;\widetilde{D_0} F)d\widetilde{W}_u
+\int^s_0 (V(r,u;\widetilde{L} F)
+2 I^{(0)}(r,u,\widetilde{D_1} F))du
\\
&=&(2\Phi(r,x_0)-1) F_0 + \hat{V}_2(r,s;F).
\end{eqnarray*}
Then we have
\begin{eqnarray*}
&&(2\Phi(r,x_0)-1) F_0 
+
\lim_{s \uparrow r} \hat{V}_2(r,s;F)
\\
&=&\lim_{s \uparrow r} V(r,s;F)
\\
&=& \widetilde{E}[(2\Phi(r-r\wedge \tau,X_{r\wedge \tau})-1)\rho_{r\wedge \tau}
F_{r\wedge \tau})|\mathcal{G}^Y_r]
\\
&=& \widetilde{E}[1_{\{\tau > r\}}(2\Phi(0,X_{r})-1)\rho_{r}
F_{r})|\mathcal{G}^Y_r]
-
\widetilde{E}[1_{\{\tau \le r\}}(2\Phi(r-\tau,0)-1)\rho_{\tau}
F_{\tau})|\mathcal{G}^Y_r]
\\
&=&
\widetilde{E}[1_{\{\tau>r\}}\rho_{r} F_r|\mathcal{G}^Y_r].
\end{eqnarray*}
So we can see that 
\begin{eqnarray*}
&&\hat{V}(r,r;F)\\
&=&\widetilde{\rho}_{r-}^{-1} e^{\int^r_0 \lambda_{x_0}(u) du} (\hat{V}_1(r,r;F) + \lambda_{x_0}(r) V(r,r;F) )
\\
&=&\widetilde{\rho}_{r-}^{-1} e^{\int^r_0 \lambda_{x_0}(u) du} 
\left\{
\hat{V}_1(r,r;F) 
+ \lambda_{x_0}(r) 
\left(
-(2\Phi(r,x_0)-1)F_0 + \widetilde{E}[1_{\{\tau>r\}}\rho_{r} F_r|\mathcal{G}^Y_r]  
\right)
\right\}
\\
&=&\hat{V}(r;f)
\end{eqnarray*}
by the first assertion of this Proposition and Proposition \ref{eq:V1}. 
\qed
\\

We now state Proposition \ref{eq:bayes}, Lemma \ref{Zakai3_1} and Proposition \ref{prop_rho} for Theorem \ref{th2}(1).

\begin{prop}\label{eq:bayes}
Let $\xi$ be a $\mathcal{B}$-measurable process. Then, we have
$$
E[\xi_{t\wedge \tau}|\mathcal{H}]
=
\frac{\widetilde{E}[\rho_{t}\xi_{t\wedge \tau}|\mathcal{H}]}{\widetilde{E}[\rho_{t}|\mathcal{H}]},\;\;\; \mathcal{H} \subset \mathcal{B}_t.
$$
\end{prop}
{\it Proof.}
For $A  \in \mathcal{H} \subset \mathcal{B}_{t}$, we have
\begin{eqnarray*}
E[\xi_{t\wedge \tau},A]
&=&
E[E[\xi_{t\wedge \tau}|\mathcal{H}],A]
=\widetilde{E}[\rho_T E[\xi_{t\wedge \tau}|\mathcal{H}],A]
\\
&=&\widetilde{E}[\widetilde{E}[\rho_T|\mathcal{B}_{t}] E[\xi_{t\wedge \tau}|\mathcal{H}],A]
=\widetilde{E}[\widetilde{E}[\rho_t|\mathcal{H}] E[\xi_{t\wedge \tau}|\mathcal{H}],A].
\end{eqnarray*}
At the same time,
$$
E[\xi_{t\wedge \tau},A]
=\widetilde{E}[\rho_T \xi_{t\wedge \tau},A]
=\widetilde{E}[\widetilde{E}[\rho_T|\mathcal{B}_{t}] \xi_{t\wedge \tau},A]
=\widetilde{E}[\rho_{t} \xi_{t\wedge \tau},A]
=\widetilde{E}[\widetilde{E}[\rho_{t} \xi_{t\wedge \tau}|\mathcal{H}],A].
$$
\qed
\begin{lemma}\label{Zakai3_1}
\begin{eqnarray*}
\widetilde{E}[\rho_t F_{t \wedge \tau}|\mathcal{F}_{t}]
&=&F_0
-\int^{t \wedge \tau}_0 \widetilde{\rho}_{r-} \hat{V}(r; F) \lambda_{x_0}(r)^{-1}  d\widetilde{M}_r
\\
&+&\int^{t \wedge \tau}_0 \widetilde{E}[\rho_{r-} (\widetilde{L} F)_r|\mathcal{F}_r] dr
+\int^{t \wedge \tau}_0 \widetilde{E}[\rho_{r-} (\widetilde{D_0} F)_r|\mathcal{F}_r]d\widetilde{W}_r.
\end{eqnarray*}
\end{lemma}
{\it Proof.}
Since $\hat{V}(r,r;F)=\hat{V}(r;F)$ by Proposition \ref{eq:V2}, we have our assertion
 by Equation (\ref{eq:rhoF}).
\qed

Let $\widetilde{\rho}_t = \widetilde{E}[\rho_t |\mathcal{F}_t]$. 
\begin{prop}\label{prop_rho}
$$
\widetilde{\rho}_t
=1
-
\int^{t}_{0} \widetilde{\rho}_{r-} \hat{V}(r;1) 
\lambda_{x_0}(r)^{-1}
d\widetilde{M}_r
+\int^{t}_{0} \widetilde{\rho}_{r-} \widetilde{E}[\beta(r,X_r,Y_r)|\mathcal{F}_r]d\widetilde{W}_r
$$
and
$$
\widetilde{\rho}^{-1}_t
=1 - \int^t_{0} \widetilde{\rho}_{r-}^{-1} \frac{\hat{V}(r;1)}{\lambda_{x_0}(r) + \hat{V}(r;1)}
 d\widetilde{M}_r
$$
$$
+\int^t_{0} \widetilde{\rho}_{r}^{-1} ( \widetilde{E}[\beta(r,X_r,Y_r)|\mathcal{F}_r]^2 
+ \frac{\hat{V}(r;1)^2 }{\lambda_{x_0}(r) + \hat{V}(r;1)}1_{\{\tau>r\}}
) dr
-\int^t_{0} \widetilde{\rho}_{r-}^{-1} \widetilde{E}[\beta(r,X_r,Y_r)|\mathcal{F}_r] d\widetilde{W}_r.
$$
\end{prop}
{\it Proof.}
Letting $F_t=1$ in Lemma \ref{Zakai3_1}, we have the first assertion.
Then we have
\begin{eqnarray*}
\widetilde{\rho}^{-1}_t 
&=& 1 - \int^t_{0} \widetilde{\rho}^{-2}_{r-} d\widetilde{\rho}_r 
+ \int^t_{0} \widetilde{\rho}_{r-}^{-3} d[\widetilde{\rho},\widetilde{\rho}]_r^c
+\sum_{0<r\le t} (\widetilde{\rho}_r^{-1} - \widetilde{\rho}_{r-}^{-1} +\widetilde{\rho}_{r-}^{-2} (\widetilde{\rho}_r - \widetilde{\rho}_{r-})
)
\\
&=&1 - \int^t_{0} \widetilde{\rho}_{r-}^{-1} \hat{V}(r;1) 
\lambda_{x_0}(r)^{-1}
d\widetilde{M}_r
-\int^t_{0} \widetilde{\rho}_{r-}^{-1}  \widetilde{E}[\beta(r,X_r,Y_r)|\mathcal{F}_r] d\widetilde{W}_r
\\
&+&
\int^t_{0} \widetilde{\rho}_{r-}^{-1}  \widetilde{E}[ \beta(r,X_r,Y_r)|\mathcal{F}_r]^2 dr
+\int^t_{0} 
\widetilde{\rho}_{r-}^{-1}
\frac{\hat{V}(r;1)^2}{\lambda_{x_0}(r) + \hat{V}(r;1)}
\lambda_{x_0}(r)^{-1}
 dN_r.
\end{eqnarray*}
Here we use the fact that
$$
\sum_{0<r\le t} (\widetilde{\rho}_r^{-1} - \widetilde{\rho}_{r-}^{-1} +\widetilde{\rho}_{r-}^{-2} (\widetilde{\rho}_r - \widetilde{\rho}_{r-}))
$$
$$
=
\sum_{0<r\le t} \frac{(\widetilde{\rho}_r - \widetilde{\rho}_{r-})^2}{\widetilde{\rho}_{r-}^2 \widetilde{\rho}_r}
=
\int^t_{0} 
\widetilde{\rho}_{r-}^{-1}
\frac{\hat{V}(r;1)^2}{\lambda_{x_0}(r)+\hat{V}(r;1)}\lambda_{x_0}(r)^{-1}dN_r.
$$
Then we have the assertion.
\qed

We give Propositions \ref{eq:Levy}, \ref{eq:FG}, and \ref{eq:limVhat} for Theorem \ref{th2}(2).

\begin{prop}\label{eq:Levy}
$$
\lambda_{a}(t) (2\Phi(t,a)-1)+\frac{\partial g}{\partial x}(t,a)=0.
$$
In particular,
$$
\hat{V}(r;F)
=
$$
$$
\widetilde{\rho}_{r-}^{-1} e^{\int^r_0 \lambda_{x_0}(u) du} 
\left(
\hat{V}_1(r,r;F) 
+
\frac{\partial g}{\partial x}(r,x_0)
 F_0 
 + \widetilde{E}[1_{\{\tau>r\}}\rho_{r} F_r|\mathcal{G}^Y_r] 
\right)
.
$$
\end{prop}
{\it Proof.}
Because of the well-known reflection principle of Brownian motion, we have $q_a(t)=P[\tau^a>t]
=1-P[\tau^a<t]=1-\frac{2}{\sqrt{2\pi t}}\int^{\infty}_{0}e^{-\frac{x^2}{2t}}dx=2\Phi(t,a)-1$.
So we have
$$
\frac{\partial}{\partial t} q_{a}(t) = 2\frac{\partial \Phi}{\partial x}(t,x)
=2\int^x_{-\infty}\frac{\partial g}{\partial t}(t,y)dy
=\frac{\partial g}{\partial x}(t,x).
$$
Since $\lambda_{a}(t) = -\frac{d}{dt}\log q_{a}(t)$, we have the assertion.
\qed
\\
\begin{prop}\label{eq:FG}
Let $Z$ be a random variable and $r>0$.Then we have
$$
\widetilde{E}[Z 1_{\{\tau>r\}}|\mathcal{G}^Y_r]1_{\{\tau>r\}}
 = 
 e^{-\int^r_0 \lambda_{x_0}(u)du} \widetilde{E}[Z|\mathcal{F}_r]1_{\{\tau>r\}}.
$$
\end{prop}
{\it Proof.}
Let $A \in \mathcal{F}_r$. Then there exists $B \in \mathcal{G}^Y_r$ such that $A \cap \{\tau>r\} = B \cap \{\tau>r\}$.
Since $\mathcal{G}^Y_r$ and $1_{\{\tau>r\}}$ are independent, we have the following.
\begin{eqnarray*}
&&\widetilde{E}[\widetilde{E}[Z 1_{\{\tau>r\}}| \mathcal{G}^Y_r]1_{\{\tau>r\}},A]\\
&=&\widetilde{E}[\widetilde{E}[Z 1_{\{\tau>r\}}| \mathcal{G}^Y_r]1_{\{\tau>r\}}1_B]
=\widetilde{E}[\widetilde{E}[Z 1_{\{\tau>r\}}1_B| \mathcal{G}^Y_r]1_{\{\tau>r\}}]
\\
&=&\widetilde{E}[\widetilde{E}[Z 1_B| \mathcal{G}^Y_r]\widetilde{E}[1_{\{\tau>r\}}]]
=\widetilde{E}[Z,A]\tilde{P}[\tau>r]
\\
&=&\widetilde{E}[ e^{-\int^r_0 \lambda_{x_0}(u)du}Z,A]
=\widetilde{E}[ e^{-\int^r_0 \lambda_{x_0}(u)du}\widetilde{E}[Z|\mathcal{F}_r],A].
\end{eqnarray*}
Then we have Assertion.
\qed
\\

\begin{prop}\label{eq:limVhat}
Let $F \in \Sigma$. Assume that there exist $C>0$ and $\alpha \in (0,1)$ such that $1_{\{|X_r| \le 1\}} 1_{\{\tau > r\}} |F_r| \le C|X_r|^{\alpha}$ for $r>0$.
Then we have $\hat{V}_1(r,r;F)=-\frac{\partial g}{\partial x}(r,x_0)F_0$ and
$$
\hat{V}(r;F)
=
\widetilde{\rho}_{r-}^{-1} e^{\int^r_0 \lambda_{x_0}(u) du}
\lambda_{x_0}(r)
 \widetilde{E}[1_{\{\tau>r\}}\rho_{r} F_r|\mathcal{G}^Y_r]
.
$$
In particular,
$$
1_{\{\tau>r\}}\hat{V}(r;F)
=
1_{\{\tau>r\}}
\widetilde{\rho}_{r-}^{-1}
\lambda_{x_0}(r)
 \widetilde{E}[1_{\{\tau>r\}}\rho_{r} F_r|\mathcal{F}_r]
.
$$
\end{prop}
{\it Proof.}
Let $1<p<\frac{2}{2-\alpha}$, $q=\frac{p}{p-1}$ and $r>s>0$.
Then we have
$$
\widetilde{E}[|I^{(1)}(r,s;\widetilde{F})|]
\le \widetilde{E}[|\frac{\partial g}{\partial x}(r-s,X_s)|1_{\{\tau > s\}} \rho_s |F_s|]
\le \widetilde{E}[1_{\{\tau>r\}}|\frac{\partial g}{\partial x}(r-s,X_s)|^p |F_s|^p]^{\frac{1}{p}}
\widetilde{E}[\rho_s^q]^{\frac{1}{q}}.
$$
Note that $x_0>0$. We have the following by Mean-Value Theorem.
$$
e^{-\frac{(x-x_0)^2}{2s}} - e^{-\frac{(x+x_0)^2}{2s}}
\le
\frac{x_0(x+x_0)}{s^2}, \;\;\; x \in (0,\infty),\;\;\; s \in (0,r).
$$
Since $\frac{2-p(2-\alpha)}{2}>0$, we have
\begin{eqnarray*}
&&\widetilde{E}[1_{\{|X_s|\le 1\}}1_{\{\tau>r\}}|\frac{\partial g}{\partial x}(r-s,X_s)|^p |F_s|^p] 
\\
&& \le
\int^{\infty}_0 (\frac{1}{\sqrt{2 \pi (r-s)}} \frac{x}{r-s} e^{-\frac{x^2}{2(r-s)}})^p
(Cx^{\alpha})^p
\frac{e^{-\frac{(x-x_0)^2}{2s}} - e^{-\frac{(x+x_0)^2}{2s}}}{\sqrt{2 \pi s}}
dx
\\
&&\le 
x_0 C^p (r-s)^{-\frac{3p}{2}}s^{-\frac{5}{2}}\int^{\infty}_0 x^{(1+\alpha)p}(x+x_0) e^{-\frac{px^2}{2(r-s)}}dx
\\
&&= 
x_0 C^p (r-s)^{\frac{2-p(2-\alpha)}{2}} s^{-\frac{5}{2}} \int^{\infty}_0 y^{(1+\alpha)p}ye^{-\frac{py^2}{2}}dy
\\
&&
+
x_0^2 C^p (r-s)^{\frac{2-p(2-\alpha)}{2}} s^{-\frac{5}{2}} \int^{\infty}_0 y^{(1+\alpha)p} e^{-\frac{py^2}{2}}dy
\\
&&
\rightarrow 0\;\;\; as \;\;\; s \uparrow r
.
\end{eqnarray*}
Let $p'>1$, $q'=\frac{p'}{p'-1}$. For $r>s>0$, we have 
\begin{eqnarray*}
&&\widetilde{E}[1_{\{|X_s| > 1\}}1_{\{\tau>r\}}|\frac{\partial g}{\partial x}(r-s,X_s)|^p |F_s|^p] 
\\
&\le& \widetilde{E}[1_{\{|X_s| > 1\}}1_{\{\tau>r\}}|\frac{\partial g}{\partial x}(r-s,X_s)|^{pp'}]^{\frac{1}{p'}}E[ |F_s|^{pq'}]^{\frac{1}{q'}} 
\\
&=& 
\left\{
\int^{\infty}_{1}  (\frac{1}{\sqrt{2 \pi (r-s)}} \frac{x}{r-s} e^{-\frac{x^2}{2(r-s)}})^{pp'}
\frac{e^{-\frac{(x-x_0)^2}{2s}} - e^{-\frac{(x+x_0)^2}{2s}}}{\sqrt{2 \pi s}}
dx
\right\}
^{\frac{1}{p'}}E[ |F_s|^{pq'}]^{\frac{1}{q'}}
\\
&\rightarrow& 0\;\;\; as \;\;\; s \uparrow r.
\end{eqnarray*}
Then we have $\hat{V}_1(r,r;F)=-\frac{\partial g}{\partial x}(r,x_0)F_0+\lim_{s \uparrow r} I^{(1)}(r,s;F)=-\frac{\partial g}{\partial x}(r,x_0)F_0$.
By Proposition \ref{eq:Levy}, we have the first assertion. We have the second assertion by Proposition \ref{eq:FG}.
\qed

We can show Theorem \ref{th2} (1) and (2) using the following proposition.
\begin{prop}\label{doubletildeMW}
$\widetilde{\widetilde{M}}_t$ is $P$-$\mathcal{F}_t$-martingale and $\widetilde{\widetilde{W}}_t$ is 
$P$-$\mathcal{B}_t$-Brownian motion.
\end{prop}
{\it Proof.}
By Proposition \ref{prop_rho},
$$
d[\widetilde{\rho}^{-1},\widetilde{M}]_t = -\widetilde{\rho}^{-1}_{t-} \frac{\hat{V}(t;1)}{\lambda_{x_0}(t) + \hat{V}(t;1)}dN_t
$$
and then we have
$$
d(\widetilde{\rho}^{-1}_t \widetilde{M}_t) 
= \widetilde{\rho}^{-1}_{t-} d\widetilde{M}_t + \widetilde{M}_{t-} d(\widetilde{\rho}^{-1})_{t} + d[\widetilde{\rho}^{-1},\widetilde{M}]_t
=\frac{\widetilde{\rho}^{-1}_{t-}\lambda_{x_0}(t)}{\lambda_{x_0}(t) + \hat{V}(t;1)}d\widetilde{\widetilde{M}}_t
+ \widetilde{M}_{t-} d(\widetilde{\rho}^{-1})_t.
$$
Since $\widetilde{\rho}^{-1}_t \widetilde{M}_t$ and $\widetilde{\rho}^{-1}_t$ are $P$-$\mathcal{F}_t$-martingale,
 we can see $\widetilde{\widetilde{M}}_t$ is also  $P$-$\mathcal{F}_t$-martingale. We can see that $\widetilde{\widetilde{W}}_t$ 
 is $P$-$\mathcal{B}_t$ -Brownian motion by the following.
 $$
d(\widetilde{\rho}^{-1}_t \widetilde{W}_t) 
= \widetilde{\rho}^{-1}_{t-} d\widetilde{W}_t + \widetilde{W}_{t} d(\widetilde{\rho}^{-1})_t + d[\widetilde{\rho}^{-1},\widetilde{W}]_t
=\widetilde{\rho}^{-1}_{t-} d\widetilde{\widetilde{W}}_t
+ \widetilde{W}_{t} d(\widetilde{\rho}^{-1})_t.
$$
\qed

Now let us prove Theorem \ref{th2}.
Let $\hat{F}_t=\widetilde{E}[\rho_{t} F_{t \wedge \tau}|\mathcal{F}_{t}]$,
 then we have the following by Lemma \ref{Zakai3_1} 
and Proposition \ref{prop_rho} .
\begin{eqnarray*}
\hat{F}_t &=& F_0 + \int^t_0 \hat{f}_0(r;F) 
d\widetilde{\widetilde{M}}_r
+ \int^t_0 \hat{f}_1(r;F) dr
+ \int^t_0 \hat{f}_2(r;F) d\widetilde{\widetilde{W}}_r,
\\
\widetilde{\rho}_t^{-1} 
&=& 1 + 
\int^t_0
 \widetilde{\rho}_{r-}^{-1} 
 \tilde{f}_0(r) 
d\widetilde{\widetilde{M}}_r
+ 
\int^t_0 
\widetilde{\rho}_{r-}^{-1} 
\tilde{f}_2(r) d\widetilde{\widetilde{W}}_r
\end{eqnarray*}
where
\begin{eqnarray*}
\hat{f}_0(r;F)
&=&
-
\frac{ \hat{V}(r;F)}{\tilde{\tilde{\lambda}}(r)-\hat{V}(r;1)}\widetilde{\rho}_{r-}
,\;\;\;
\hat{f}_2(r;F)=\widetilde{E}[\rho_{r-}(\widetilde{D_0} F)_r|\mathcal{F}_r] ,
\\
\hat{f}_1(r;F)&=&
\widetilde{E}[\rho_{r-}(\widetilde{L} F)_r|\mathcal{F}_r]
\\
&+&1_{\{\tau>r\}}
\frac{\hat{V}(r;1)\hat{V}(r;F) }
{\tilde{\tilde{\lambda}}(r) - \hat{V}(r;1)}\widetilde{\rho}_{r-}
+ E[\beta(r,X_{r \wedge \tau},Y_r)|\mathcal{F}_r]\widetilde{E}[\rho_{r-}(\widetilde{D_0} F)_r|\mathcal{F}_r],
\\
\tilde{f}_0(r) &=& -\hat{V}(r;1)\tilde{\tilde{\lambda}}(r)^{-1},\;\;\;
\tilde{f}_2(r) = -E[\beta(r,X_r,Y_r)|\mathcal{F}_r].
\end{eqnarray*}
Note that $d\widetilde{M}_t = d\widetilde{\widetilde{M}}_t + (1 - N_{t-})\hat{V}(t;1)dt$
 and $d\widetilde{W}_t = d\widetilde{\widetilde{W}}_t + E[\beta(t,X_{t \wedge \tau},Y_t)|\mathcal{F}_t]dt$.
Then $E[F_{t \wedge \tau}|\mathcal{F}_{t}]=\widetilde{\rho}_t^{-1} \hat{F}_t$. Let $\bar{F}_t = E[F_{t \wedge \tau}|\mathcal{F}_{t}]$ and we have the following.
\begin{eqnarray*}
&&\widetilde{\rho}_t^{-1} \hat{F}_t 
\\
&=& F_0 
+
\int^{t\wedge \tau}_{0} \widetilde{\rho}_{r-}^{-1} d\hat{F}_r 
+
\int^{t\wedge \tau}_{0} \hat{F}_{r-} d \widetilde{\rho}_r^{-1}
+
[\hat{F},\widetilde{\rho}^{-1}]_{t\wedge \tau}
\\
&=& F_0 
+ 
\int^{t\wedge \tau}_{0} \widetilde{\rho}_{r-}^{-1}
\left(
\hat{f}_0(r;F) +  \tilde{f}_0(r)\hat{F}_{r-}
\right)
d\widetilde{\widetilde{M}}_r
+
\int^{t\wedge \tau}_{0} \widetilde{\rho}_{r-}^{-1}
\left(
\hat{f}_1(r;F) + \tilde{f}_2(r)\hat{f}_2(r;F)
\right)
dr
\\
&+&
\int^{t\wedge \tau}_{0} 
\widetilde{\rho}_{r-}^{-1}
\left(
\hat{f}_2(r;F) +  \tilde{f}_2(r)\hat{F}_{r-}
\right)
d\widetilde{\widetilde{W}}_r
+
\sum_{0< r \le {t\wedge \tau}} (\widetilde{\rho}_r^{-1} - \widetilde{\rho}_{r-}^{-1} )(\hat{F}_r -\hat{F}_{r-} )
\\
&=& F_0 
+ \int^{t\wedge \tau}_{0} \widetilde{\rho}_{r-}^{-1}
\left(
\hat{f}_0(r;F) + \tilde{f}_0(r)\hat{F}_{r-}
+
\tilde{f}_0(r)\hat{f}_0(r;F)
\right)
d\widetilde{\widetilde{M}}_r
\\
&+&
\int^{t\wedge \tau}_{0} \widetilde{\rho}_{r-}^{-1}
\left(
\hat{f}_1(r;F) +\tilde{f}_2(r)\hat{f}_2(r;F)
+ 1_{\{\tau>r\}}
\tilde{f}_0(r)
\hat{f}_0(r;F)
\tilde{\tilde{\lambda}}(r)
\right)
dr
\\
&+& 
\int^{t\wedge \tau}_{0} \widetilde{\rho}_{r-}^{-1}
\left(
\hat{f}_2(r;F) +  \tilde{f}_2(r)\hat{F}_{r-}
\right)
d\widetilde{\widetilde{W}}_r.
\end{eqnarray*}
Here we note that
$$
\sum_{0< r \le t} (\widetilde{\rho}_r^{-1} - \widetilde{\rho}_{r-}^{-1} )(\hat{F}_r -\hat{F}_{r-} )
=
\int^t_0
\widetilde{\rho}_{r-}^{-1} 
\tilde{f}_0(r)
\hat{f}_0(r;F)dN_r.
$$
Then we have Theorem \ref{th2}(1) as the following.
\begin{eqnarray*}
E[F_{t \wedge \tau}|\mathcal{F}_{t}]
&=&F_0
-\int^t_{0}  
1_{\{\tau>r\}}
\left(
\hat{V}(r; F) + \hat{V}(r;1)\bar{F}_{r-}
\right)
\tilde{\tilde{\lambda}}(r)^{-1}
 d\widetilde{\widetilde{M}}_r
\\
&+&\int^t_{0}
 1_{\{\tau>r\}}
 E[1_{\{\tau>r\}} (\widetilde{L} F)_r|\mathcal{F}_r]
dr
\\
&+&
\int^t_{0} 
\left(
E[(\widetilde{D_0} F)_r|\mathcal{F}_r]
 - 
E[\beta(r,X_r,Y_r)|\mathcal{F}_r]
\bar{F}_{r-}
\right) 
d\widetilde{\widetilde{W}}_r.
\end{eqnarray*}
If there exist $C>0$ and $\alpha \in (0,1)$  such that $1_{\{|X_t| \le 1\}}1_{\{\tau > t\}}|F_t| \le C|X_t|^{\alpha}$ for $t>0$,
we have
\begin{eqnarray*}
1_{\{\tau>r\}}\hat{V}(r;F)
&=&
1_{\{\tau>r\}}
\widetilde{\rho}_{r-}^{-1}
\lambda_{x_0}(r)
 \widetilde{E}[1_{\{\tau>r\}}\rho_{r} F_r|\mathcal{F}_r]
\\
&=&1_{\{\tau>r\}}
\widetilde{\rho}_{r-}^{-1}
\lambda_{x_0}(r)
e^{\int^r_0 \lambda_{x_0} (u) du}
 \widetilde{E}[1_{\{\tau>r\}}\rho_{r} F_r|\mathcal{G}_r]
 \end{eqnarray*}
by Proposition \ref{eq:FG} and Proposition \ref{eq:limVhat}. Then we have
\begin{eqnarray*}
\tilde{f}_0(r;F)
&=&
-1_{\{\tau>r\}}
\left(
\hat{V}(r; F) + \hat{V}(r;1)\bar{F}_{r-}
\right)
\tilde{\tilde{\lambda}}(r)^{-1}
\\
&=&
-1_{\{\tau>r\}}
\frac{
\lambda_{x_0}(r)\widetilde{\rho}_{r-}^{-1}
e^{\int^r_0 \lambda_{x_0} (u) du}
 \widetilde{E}[1_{\{\tau>r\}}\rho_{r} F_r|\mathcal{G}_r]
+
 \hat{V}(r;1)\bar{F}_{r-}
}
{\lambda_{x_0}(r)+\hat{V}(r;1)}
\\
&=&
-1_{\{\tau>r\}}
\bar{F}_{r-},
\end{eqnarray*}
which gives Theorem \ref{th2}(2).
%
\vskip24pt
\small
%
\centerline{\bf References}
\vskip12pt
\begin{enumerate}
\renewcommand{\labelenumi}{[\arabic{enumi}]}
\renewcommand{\makelabel}{\rm}
\setcounter{enumi}{0}
\setlength{\itemsep}{-3pt}
\setlength{\parsep}{0cm}

\bibitem{Rutkowski} Bielecki, Tomasz R., and Marek Rutkowski, 2002 ``Credit risk: modeling, valuation and hedging." Springer, 2002.

\bibitem{DuffieLando} Duffie, Darrell, and David Lando, 2001 ``Term structures of credit spreads with incomplete accounting information." Econometrica 69.3: 633-664.

\bibitem{JarrowProtterDeniz} Jarrow, Robert A., Philip Protter, and A. Deniz Sezer, 2007. ``International Journal of Theoretical and Applied Finance." Finance and Stochastics 11.2: 195-212.

\bibitem{Jeanblanc} Jeanblanc, Monique, and Stoyan Valchev, 2005 ``Partial information and hazard process." International Journal of Theoretical and Applied Finance 8.06: 807-838.

\bibitem{Nakagawa} Nakagawa, Hidetoshi, 2001 ``A Filtering Model on Default Risk." J. Math. Sci. Univ. Tokyo 8: 107-142.

\end{enumerate}

\vskip12pt

\parindent7cm


\label{finalpage}
\end{document}